\documentclass[12pt,vatola]{article}

\usepackage{graphicx}

\def\c{\centerline}

\def\no{\noindent}

\begin{document}

\c{\Large\bf Smarandache Multi-Space Theory(III) }\vskip 5mm

\hskip 50mm {\it -Map geometries and pseudo-plane
geometries}\vskip 10mm

\c{Linfan Mao}\vskip 3mm

\c{\small Academy of Mathematics and System Sciences}

\c{\small Chinese Academy of Sciences, Beijing 100080}

\c{\small maolinfan@163.com}

\vskip 10mm

\begin{minipage}{130mm}

\no{\bf Abstract.} {\small A Smarandache multi-space is a union of
$n$ different spaces equipped with some different structures for
an integer $n\geq 2$, which can be both used for discrete or
connected spaces, particularly for geometries and spacetimes in
theoretical physics. This monograph concentrates on characterizing
various multi-spaces including three parts altogether. The first
part is on {\it algebraic multi-spaces with structures}, such as
those of multi-groups, multi-rings, multi-vector spaces,
multi-metric spaces, multi-operation systems and multi-manifolds,
also multi-voltage graphs, multi-embedding of a graph in an
$n$-manifold,$\cdots$, etc.. The second discusses {\it Smarandache
geometries}, including those of map geometries, planar map
geometries and pseudo-plane geometries, in which the {\it Finsler
geometry}, particularly the {\it Riemann geometry} appears as a
special case of these Smarandache geometries. The third part of
this book considers the {\it applications of multi-spaces to
theoretical physics}, including the relativity theory, the
M-theory and the cosmology. Multi-space models for $p$-branes and
cosmos are constructed and some questions in cosmology are
clarified by multi-spaces. The first two parts are relative
independence for reading and in each part open problems are
included for further research of interested readers.}

\vskip 3mm \no{\bf Key words:} {\small  Smarandache geometries,
map geometries, planar map geometries, pseudo-plane geometries,
Finsler geometries.}

 \vskip 3mm \no{{\bf
Classification:} AMS(2000) 03C05,05C15,51D20,51H20,51P05,83C05,
83E50}
\end{minipage}

\newpage

\no{\large\bf Contents}\vskip 8mm

\no $3.$ Map geometries\dotfill 3\vskip 3mm

\no $\S 3.1$ \ Smarandache Geometries\dotfill 3\vskip 2mm

\no $3.1.1$ What are lost in classical mathematics\dotfill 3

\no $3.1.2$ Smarandache geometries\dotfill 5

\no $3.1.3$ Smarandache manifolds\dotfill 8\vskip 2mm

\no $\S 3.2$ \ Map Geometries without Boundary\dotfill 9\vskip 2mm

\no $\S 3.3$ \ Map Geometries with Boundary\dotfill 18\vskip 2mm

\no $\S 3.4$ \ The Enumeration of Map Geometries\dotfill 21\vskip
2mm

\no $\S 3.5$ \ Remarks and Open Problems\dotfill 24

\vskip 5mm

\no $4.$ Planar map geometries\dotfill 26\vskip 3mm

\no $\S 4.1$ \ Points in a Planar Map Geometry\dotfill 26\vskip
2mm

\no $\S 4.2$ \ Lines in a Planar Map Geometry \dotfill 29\vskip
2mm

\no $4.2.1$ Lines in a planar map geometry \dotfill 29

\no $4.2.2$ Curvature of an $m$-line \dotfill 32\vskip 2mm

\no $\S 4.3$ \ Polygons in a Planar Map Geometry \dotfill 35\vskip
2mm

\no $4.3.1$ Existence \dotfill 35

\no $4.3.2$ Sum of internal angles \dotfill 37

\no $4.3.3$ Area of a polygon \dotfill 39\vskip 2mm

\no $\S 4.4$ \ Circles in a Planar Map Geometry\dotfill 41\vskip
2mm

\no $\S 4.5$ \ Line Bundles in a Planar Map Geometry\dotfill
43\vskip 2mm

\no $4.5.1$ A condition for parallel bundles\dotfill 44

\no $4.5.2$ Linear conditions and combinatorial realization for
parallel bundles\dotfill 49\vskip 2mm

\no $\S 4.6$ \ Examples of Planar Map Geometries\dotfill 52\vskip
2mm

\no $\S 4.7$ \ Remarks and Open Problems\dotfill 55\vskip 3mm

\vskip 5mm

\no $5.$ Pseudo-Plane geometries\dotfill 57\vskip 3mm

\no $\S 5.1$ \ Pseudo-Planes\dotfill 57\vskip 2mm

\no $\S 5.2$ \ Integral Curves\dotfill 62\vskip 2mm

\no $\S 5.3$ \ Stability of a Differential Equation\dotfill
67\vskip 2mm

\no $\S 5.4$ \ Remarks and Open Problems\dotfill 71\vskip 3mm

\newpage

\no{\large\bf $3.$ Map geometries}\vskip 25mm

\no As a kind of multi-metric spaces, Smarandache geometries were
introduced by Smarandache in $[86]$ and investigated by many
mathematicians. These geometries are related with the {\it Euclid
geometry}, the {\it Lobachevshy-Bolyai-Gauss geometry} and the
{\it Riemann geometry}, also related with {\it relativity theory}
and {\it parallel universes} (see $[56],[35]-[36],[38]$ and
$[77]-[78]$ for details). As a generalization of Smarandache
manifolds of dimension $2$, {\it Map geometries} were introduced
in $[55],[57]$ and $[62]$, which can be also seen as a realization
of Smarandache geometries on surfaces or Smarandache geometries on
maps.

\vskip 5mm

\no{\bf \S $3.1$ \ Smarandache Geometries}

\vskip 4mm

\no{\bf $3.1.1.$ What are lost in classical mathematics?}

\vskip 3mm

\no As we known, mathematics is a powerful tool of sciences for
its unity and neatness, without any shade of mankind. On the other
hand, it is also a kind of aesthetics deep down in one's mind.
There is a famous proverb says that {\it only the beautiful things
can be handed down to today}, which is also true for the
mathematics.

Here, the terms {\it unity} and {\it neatness} are relative and
local, maybe also have various conditions. For obtaining a good
result, many unimportant matters are abandoned in the research
process. Whether are those matters still unimportant in another
time? It is not true. That is why we need to think a queer
question: {\it what are lost in the classical mathematics?}

For example, a compact surface is topological equivalent to a
polygon with even number of edges by identifying each pairs of
edges along its a given direction ($[68],[92]$). If label each
pair of edges by a letter $e,e\in {\mathcal E}$, a surface $S$ is
also identified to a cyclic permutation such that each edge
$e,e\in {\mathcal E}$ just appears two times in $S$, one is $e$
and another is $e^{-1}$ (orientable) or $e$ (non-orientable). Let
$a,b,c,\cdots$ denote letters in ${\mathcal E}$ and $A,B,C,\cdots$
the sections of successive letters in a linear order on a surface
$S$ (or a string of letters on $S$). Then, an orientable surface
can be represented by

$$S=(\cdots , A,a,B,a^{-1},C,\cdots),$$

\no{where£¬$a\in {\mathcal E}$ and $A,B,C$ denote strings of
letter. Three elementary transformations are defined as follows:}

\vskip 2mm

$(O_1)\quad\quad (A,a,a^{-1},B)\Leftrightarrow (A,B);$

\vskip 2mm

$(O_2)\quad\quad (i)\quad (A,a,b,B,b^{-1},a^{-1})\Leftrightarrow
(A,c,B,c^{-1}) ;$

\quad\quad\quad\quad $(ii)\quad (A,a,b,B,a,b)\Leftrightarrow
(A,c,B,c); $

\vskip 2mm

$(O_3)\quad\quad (i)\quad (A,a,B,C,a^{-1},D)\Leftrightarrow
(B,a,A,D,a^{-1},C);$

\quad\quad\quad\quad $(ii)\quad (A,a,B,C,a,D)\Leftrightarrow
(B,a,A,C^{-1},a,D^{-1}).$ \vskip 2mm

\no If a surface $S_0$ can be obtained by these elementary
transformations $O_1$-$O_3$ from a surface $S$, it is said that
$S$ is {\it elementary equivalent} with $S_0$, denoted by
$S\sim_{El}S_0$.

We have known the following formulae from $[43]$:\vskip 2mm

$(i)\quad
(A,a,B,b,C,a^{-1},D,b^{-1},E)\sim_{El}(A,D,C,B,E,a,b,a^{-1},b^{-1});$

$(ii)\quad (A,c,B,c)\sim_{El} (A,B^{-1},C,c,c);$

$(iii)\quad (A,c,c,a,b,a^{-1},b^{-1})\sim_{El} (A,c,c,a,a,b,b).$
\vskip 2mm

\no Then we can get the classification theorem of compact surfaces
as follows($[68]$):\vskip 2mm

{\it Any compact surface is homeomorphic to one of the following
standard surfaces:

($P_0$) \ The sphere: $aa^{-1}$;

($P_n$) \ The connected sum of $n,n\geq 1$, tori:

$$a_1b_1a_1^{-1}b_1^{-1}a_2b_2a_2^{-1}b_2^{-1}\cdots a_nb_na_n^{-1}b_n^{-1};$$

($Q_n$) \ The connected sum of $n,n\geq 1$, projective planes:

$$a_1a_1a_2a_2\cdots a_na_n.$$}\vskip 2mm

As we have discussed in Chapter $2$, a combinatorial map is just a
kind of decomposition of a surface. Notice that all the standard
surfaces are one face map underlying an one vertex graph, i.e., a
bouquet $B_n$ with $n\geq 1$. By a combinatorial view, {\it a
combinatorial map is nothing but a surface}. This assertion is
needed clarifying. For example, let us see the left graph $\Pi_4$
in Fig. $3.1$, which is a tetrahedron.

\includegraphics[bb=10 10 200 140]{sgm37.eps}

\vskip 2mm

\c{\bf Fig.$3.1$}

\vskip 3mm

\no Whether can we say $\prod_4$ is a sphere? Certainly NOT. Since
any point $u$ on a sphere has a neighborhood $N(u)$ homeomorphic
to an open disc, thereby all angles incident with the point $1$
must be $120^{\circ}$ degree on a sphere. But in $\Pi_4$, those
are only $60^{\circ}$ degree. For making them same in a
topological sense, i.e., homeomorphism, we must blow up the
$\Pi_4$ and make it become a sphere. This physical processing is
shown in the Fig.$3.1$. Whence, for getting the classification
theorem of compact surfaces, we lose the {\it angle,area,
volume,distance,curvature},$\cdots$, etc, which are also lost in
combinatorial maps.

By a geometrical view, {\it Klein Erlanger Program} says that {\it
any geometry is nothing but find invariants under a transformation
group of this geometry}. This is essentially the group action idea
and widely used in mathematics today. Surveying topics appearing
in publications for combinatorial maps, we know the following
problems are applications of {\it Klein Erlanger Program}:

\vskip 3mm

 ($i$) \ {\it to determine isomorphism maps or rooted maps;}

($ii$) \ {\it to determine equivalent embeddings of a graph;}

($iii$) \ {\it to determine an embedding whether exists or not;}

($iv$) \ {\it to enumerate maps or rooted maps on a surface;}

($v$) \ {\it to enumerate embeddings of a graph on a surface;}

($vi$) \ $\cdots$, etc. \vskip 2mm

All the problems are extensively investigated by researches in the
last century and papers related those problems are still
frequently appearing in journals today. Then,\vskip 3mm

{\it what are their importance to classical mathematics?}\vskip
2mm

\no and\vskip 3mm

{\it what are their contributions to sciences?} \vskip 2mm

Today, we have found that combinatorial maps can contribute an
underlying frame for applying mathematics to sciences, i.e.,
through by map geometries or by graphs in spaces.

\vskip 4mm

\no{\bf $3.1.2.$ Smarandache geometries}

\vskip 3mm

\no{\it Smarandache geometries} were proposed by Smarandache in
$[86]$ which are generalization of classical geometries, i.e.,
these {\it Euclid, Lobachevshy-Bolyai-Gauss} and {\it Riemann
geometries} may be united altogether in a same space, by some
Smarandache geometries. These last geometries can be either
partially Euclidean and partially Non-Euclidean, or Non-Euclidean.
Smarandache geometries are also connected with the {\it Relativity
Theory} because they include Riemann geometry in a subspace and
with the {\it Parallel Universes} because they combine separate
spaces into one space too. For a detail illustration, we need to
consider classical geometries first.

As we known, the axiom system of an {\it Euclid geometry} is in
the following:

\vskip 3mm

(A1) \ {\it there is a straight line between any two points.}

(A2) \ {\it  a finite straight line can produce a infinite
straight line continuously.}

(A3) \ {\it any point and a distance can describe a circle.}

(A4) \ {\it all right angles are equal to one another.}

(A5) \ {\it if a straight line falling on two straight lines make
the interior angles on the same side less than two right angles,
then the two straight lines, if produced indefinitely, meet on
that side on which are the angles less than the two right angles.}

\vskip 2mm

The axiom (A5) can be also replaced by: \vskip 3mm

(A5') \ {\it  given a line and a point exterior this line, there
is one line parallel to this line.}

\vskip 2mm

The {\it Lobachevshy-Bolyai-Gauss geometry}, also called {\it
hyperbolic geometry}, is a geometry with axioms $(A1)-(A4)$ and
the following axiom $(L5)$:\vskip 3mm

(L5) \ {\it there are infinitely many lines parallel to a given
line passing through an exterior point.}\vskip 2mm

The {\it Riemann geometry}, also called {\it elliptic geometry},
is a geometry with axioms $(A1)-(A4)$ and the following axiom
$(R5)$:\vskip 3mm

{\it there is no parallel to a given line passing through an
exterior point.} \vskip 2mm

By a thought of anti-mathematics: {\it not in a nihilistic way,
but in a positive one, i.e., banish the old concepts by some new
ones: their opposites}, Smarandache introduced thse {\it
paradoxist geometry, non-geometry, counter-projective geometry}
and {\it anti-geometry} in $[86]$ by contradicts axioms
$(A1)-(A5)$ in an Euclid geometry.

\vskip 4mm

\no{\bf Paradoxist geometry}

\vskip 3mm

\no In this geometry, its axioms consist of $(A1)-(A4)$ and one of
the following as the axiom $(P5)$: \vskip 3mm

($i$) \ {\it there are at least a straight line and a point
exterior to it in this space for which any line that passes
through the point intersect the initial line.}

($ii$) \ {\it there are at least a straight line and a point
exterior to it in this space for which only one line passes
through the point and does not intersect the initial line.}

($iii$) \ {\it there are at least a straight line and a point
exterior to it in this space for which only a finite number of
lines $l_1,l_2,\cdots , l_k, k\geq 2$ pass through the point and
do not intersect the initial line.}

($iv$) \ {\it there are at least a straight line and a point
exterior to it in this space for which an infinite number of lines
pass through the point (but not all of them) and do not intersect
the initial line.}

($v$) \ {\it there are at least a straight line and a point
exterior to it in this space for which any line that passes
through the point and does not intersect the initial line.}\vskip
2mm

\vskip 4mm

\no{\bf Non-Geometry}

\vskip 3mm

\no The non-geometry is a geometry by denial some axioms of
$(A1)-(A5)$, such as:

\vskip 3mm

($A1^-$) \ {\it It is not always possible to draw a line from an
arbitrary point to another arbitrary point.}

($A2^-$) \ {\it It is not always possible to extend by continuity
a finite line to an infinite line.}

($A3^-$) \ {\it It is not always possible to draw a circle from an
arbitrary point and of an arbitrary interval.}

($A4^-$) \ {\it not all the right angles are congruent.}

($A5^-$) \ {\it if a line, cutting two other lines, forms the
interior angles of the same side of it strictly less than two
right angle, then not always the two lines extended towards
infinite cut each other in the side where the angles are strictly
less than two right angle}.

\vskip 4mm

\no{\bf Counter-Projective geometry}

\vskip 3mm

\no Denoted by $P$ the point set, $L$ the line set and $R$ a
relation included in $P\times L$. A counter-projective geometry is
a geometry with these counter-axioms $(C_1)-(C_3)$:\vskip 3mm

($C1$) \ {\it there exist: either at least two lines, or  no line,
that contains two given distinct points.}

($C2$) \ {\it let $p_1,p_2,p_3$ be three non-collinear points, and
$q_1,q_2$ two distinct points. Suppose that $\{p_1.q_1,p_3\}$ and
$\{p_2,q_2,p_3\}$ are collinear triples. Then the line containing
$p_1,p_2$ and the line containing $q_1,q_2$ do not intersect.}

($C3$) \ {\it every line contains at most two distinct points. }

\vskip 4mm

\no{\bf Anti-Geometry}

\vskip 3mm

\no A geometry by denial some axioms of the Hilbert's $21$ axioms
of Euclidean geometry. As shown in $[38]$, there are at least
$2^{21}-1$ anti-geometries.

In general, Smarandache geometries are defined as follows.

\vskip 4mm

\no{\bf Definition $3.1.1$} {\it An axiom is said to be
Smarandachely denied if the axiom behaves in at least two
different ways within the same space, i.e., validated and
invalided, or only invalided but in multiple distinct ways.

A Smarandache geometry is a geometry which has at least one
Smarandachely denied axiom($1969$).}

\vskip 3mm

In a Smarandache geometries, points, lines, planes, spaces,
triangles, $\cdots$, etc are called $s$-points, $s$-lines,
$s$-planes, $s$-spaces, $s$-triangles, $\cdots$, respectively in
order to distinguish them from classical geometries. An example of
Smarandache geometries in the classical geometrical sense is in
the following.

\vskip 3mm

\no{\bf Example $3.1.1$} \ Let us consider an Euclidean plane
${\bf R}^2$ and three non-collinear points $A,B$ and $C$. Define
$s$-points as all usual Euclidean points on ${\bf R}^2$ and
$s$-lines any Euclidean line that passes through one and only one
of points $A,B$ and $C$. Then this geometry is a Smarandache
geometry because two axioms are Smarandachely denied comparing
with an Euclid geometry:

\vskip 2mm

($i$) \ The axiom (A5) that through a point exterior to a given
line there is only one parallel passing through it is now replaced
by two statements: {\it one parallel}, and {\it no parallel}. Let
$L$ be an $s$-line passes through $C$ and is parallel in the
euclidean sense to $AB$. Notice that through any $s$-point not
lying on $AB$ there is one $s$-line parallel to $L$ and through
any other $s$-point lying on $AB$ there is no $s$-lines parallel
to $L$ such as those shown in Fig.$3.2(a)$.

\includegraphics[bb=10 10 200 120]{sgm38.eps}

\vskip 2mm

\c{\bf Fig.$3.2$} \vskip 3mm

($ii$) \ The axiom that through any two distinct points there
exist one line passing through them is now replaced by; {\it one
$s$-line}, and {\it no $s$-line}. Notice that through any two
distinct $s$-points $D,E$ collinear with one of $A,B$ and $C$,
there is one $s$-line passing through them and through any two
distinct $s$-points $F,G$ lying on $AB$ or non-collinear with one
of $A,B$ and $C$, there is no $s$-line passing through them such
as those shown in Fig.$3.2(b)$.

\vskip 4mm

\no{\bf $3.1.3.$ Smarandache manifolds}

\vskip 3mm

\no Generally, a {\it Smarandache manifold} is an $n$-dimensional
manifold that support a Smarandache geometry. For $n=2$, a nice
model for Smarandache geometries called {\it $s$-manifolds} was
found by Iseri in $[35][36]$, which is defined as follows:\vskip
3mm

{\it An $s$-manifold is any collection ${\mathcal C}(T,n)$ of
these equilateral triangular disks $T_i, 1\leq i\leq n$ satisfying
the following conditions:}

$(i)$ \ {\it each edge $e$ is the identification of at most two
edges $e_i,e_j$ in two distinct triangular disks $T_i,T_j, 1\leq
i,j\leq n$ and $i\not= j$;}

$(ii)$ \ {\it each vertex $v$ is the identification of one vertex
in each of five, six or seven distinct triangular disks.}\vskip
2mm

The vertices are classified by the number of the disks around
them. A vertex around five, six or seven triangular disks is
called an {\it elliptic vertex}, an {\it euclidean vertex} or a
{\it hyperbolic vertex}, respectively.

In a plane, an elliptic vertex $O$, an euclidean vertex $P$ and a
hyperbolic vertex $Q$ and an $s$-line $L_1$, $L_2$ or $L_3$ passes
through points $O,P$ or $Q$ are shown in Fig.$3.3(a),(b),(c)$,
respectively.

\includegraphics[bb=0 10 400 120]{sgm39.eps}

\vskip 2mm

\c{\bf Fig.$3.3$} \vskip 3mm

Smarandache paradoxist geometries and non-geometries can be
realized by $s$-manifolds, but other Smarandache geometries can be
only partly realized by this kind of manifolds. Readers are
inferred to Iseri's book $[35]$ for those geometries.

An $s$-manifold is called closed if each edge is shared exactly by
two triangular disks. An elementary classification for closed
$s$-manifolds by triangulation were introduced in $[56]$. They are
classified into $7$ classes. Each of those classes is defined in
the following.

\vskip 3mm

{\bf Classical Type}:\vskip 2mm

$(1)$ \ $\Delta_1=\{5-regular \ triangular \ maps\}$ ({\it
elliptic});

$(2)$ \ $\Delta_2=\{6-regular  \ triangular \ maps\}$({\it
euclidean});

$(3)$ \ $\Delta_3=\{7-regular  \ triangular \ maps\}$({\it
hyperbolic}).\vskip 3mm

{\bf Smarandache Type}:\vskip 2mm

$(4)$ \ $\Delta_4=\{triangular \ maps \ with \ vertex \ valency \
5 \ and \ 6\}$ ({\it euclid-elliptic});

$(5)$ \ $\Delta_5=\{triangular \ maps \ with \ vertex \ valency \
5 \ and \ 7\}$ ({\it elliptic-hyperbolic});

$(6)$ \ $\Delta_6=\{triangular \ maps \ with \ vertex \ valency \
6 \ and \ 7\}$ ({\it euclid-hyperbolic});

$(7)$ \ $\Delta_7=\{triangular \ maps \ with \ vertex \ valency \
5 , 6 \ and \ 7\}$ ({\it mixed}).\vskip 2mm

It is proved in $[56]$ that $|\Delta_1|=2$, $|\Delta_5|\geq 2$ and
$|\Delta_i|, i=2,3,4,6,7$ are infinite. Iseri proposed a question
in $[35]$: {\it Do the other closed $2$-manifolds correspond to
$s$-manifolds with only hyperbolic vertices?} Since there are
infinite Hurwitz maps, i.e., $|\Delta_3|$ is infinite, the answer
is affirmative.

\vskip 5mm

\no{\bf \S $3.2$ \ Map Geometries without Boundary}

\vskip 4mm

\no A combinatorial map can be also used to construct new models
for Smarandache geometries. By a geometrical view, these models
are generalizations of Isier's model for Smarandache geometries.
For a given map on a locally orientable surface, map geometries
without boundary are defined in the following definition.

\vskip 4mm

\no{\bf Definition $3.2.1$} \ {\it For a combinatorial map $M$
with each vertex valency$\geq 3$, associates a real number $\mu
(u), 0 \ < \mu (u) \ < \ \frac{4\pi}{\rho_M(u)}$, to each vertex
$u, u\in V(M)$. Call $(M,\mu)$ a map geometry without boundary,
$\mu (u)$ an angle factor of the vertex $u$ and orientablle or
non-orientable if $M$ is orientable or not.}

\vskip 3mm

The realization for vertices $u,v,w\in V(M)$ in a space ${\bf
R}^3$ is shown in Fig.$3.4$, where $\rho_M(u)\mu (u) <  2\pi$ for
the vertex $u$, $\rho_M(v)\mu (v) = 2\pi$ for the vertex $v$ and
$\rho_M(w)\mu (w)
> 2\pi$ for the vertex $w$, respectively.

\vskip 2mm

\includegraphics[bb=30 20 100 100]{sgm40.eps}

\vskip 2mm

\c{$\rho_M(u)\mu (u) \ < \ 2\pi$\hskip 18mm $\rho_M(u)\mu (u)
=2\pi$\hskip 20mm $\rho_M(u)\mu (u) \ > \ 2\pi$}\vskip 2mm

\c{\bf Fig.$3.4$}

As we have pointed out in Section $3.1$, this kind of realization
is not a surface, but it is homeomorphic to a locally orientable
surface by a view of topological equivalence. Similar to
$s$-manifolds, we also classify points in a map geometry $(M,\mu)$
without boundary into {\it elliptic points}, {\it euclidean
points} and {\it hyperbolic points}, defined in the next
definition.

\vskip 4mm

\no{\bf Definition $3.2.2$} \ {\it A point $u$ in a map geometry
$(M,\mu)$ is said to be elliptic, euclidean or hyperbolic if
$\rho_M (u)\mu (u)  <  2\pi$, $\rho_M (u)\mu (u) =2\pi$ or $\rho_M
(u)\mu (u)  >  2\pi$.}\vskip 3mm

Then we get the following results.

\vskip 4mm

\no{\bf Theorem $3.2.1$} \ {\it Let $M$  be a map with $\forall
u\in V(M), \rho_M(u)\geq 3$. Then for $\forall u\in V(M)$, there
is a map geometry $(M,\mu)$ without boundary such that $u$ is
elliptic, euclidean or hyperbolic.}

\vskip 3mm

{\it Proof} \ Since $\rho_M(u)\geq 3$, we can choose an angle
factor $\mu (u)$ such that $\mu (u)\rho_M(u)  <  2\pi$, $\mu
(u)\rho_M(u)= 2\pi$ or $\mu (u)\rho_M(u)  >  2\pi$. Notice that

$$0 \ < \frac{2\pi}{\rho_M (u)}  <  \frac{4\pi}{\rho_M(u)}.$$

\no Thereby we can always choose $\mu (u)$ satisfying that $ 0 \ <
\mu (u) \ < \frac{4\pi}{\rho_M(u)}. \quad\quad \natural$

\vskip 4mm

\no{\bf Theorem $3.2.2$} \ {\it Let $M$  be a map of order$\geq 3$
and $\forall u\in V(M), \rho_M (u)\geq 3$. Then there exists a map
geometry $(M,\mu)$ without boundary in which elliptic, euclidean
and hyperbolic points appear simultaneously.}

\vskip 3mm

{\it Proof} \ According to Theorem $3.2.1$,  we can always choose
an angle factor $\mu$ such that a vertex $u, u\in V(M)$ to be
elliptic, or euclidean, or hyperbolic. Since $|V(M)|\geq 3$, we
can even choose the angle factor $\mu$ such that any two different
vertices $v,w\in V(M)\backslash \{u\}$ to be elliptic, or
euclidean, or hyperbolic as we wish. Then the map geometry
$(M,\mu)$ makes the assertion hold.\quad\quad $\natural$

A {\it geodesic} in a manifold is a curve as straight as possible.
Applying conceptions such as angles and straight lines in an
Euclid geometry, we define $m$-lines and $m$-points in a map
geometry in the next definition.

\vskip 4mm

\no{\bf Definition $3.2.3$} \ {\it Let $(M,\mu)$ be a map geometry
without boundary and let $S(M)$ be the locally orientable surface
represented by a plane polygon on which $M$ is embedded. A point
$P$ on $S(M)$ is called an $m$-point. A line $L$ on $S(M)$ is
called an $m$-line if it is straight in each face of $M$ and each
angle on $L$ has measure $\frac{\rho_M(v)\mu(v)}{2}$ when it
passes through a vertex $v$ on $M$.}

\vskip 3mm

Two examples for $m$-lines on the torus are shown in the
Fig.$3.5(a)$ and $(b)$, where $M=M(B_2)$, $\mu(u)=\frac{\pi}{2}$
for the vertex $u$ in $(a)$ and

$$\mu(u)=\frac{135-\arctan(2)}{360}\pi$$

\no for the vertex $u$ in $(b)$, i.e., $u$ is euclidean in $(a)$
but elliptic in $(b)$. Notice that in $(b)$, the $m$-line $L_2$ is
self-intersected.

\vskip 3mm

\includegraphics[bb=10 10 100 140]{sgm41.eps}

\vskip 2mm

\c{\bf Fig.$3.5$}\vskip 2mm

If an $m$-line passes through an elliptic point or a hyperbolic
point $u$, it must has an angle $\frac{\mu (u)\rho_M(u)}{2}$ with
the entering line, not $180^{\circ}$ which are explained in
Fig.$3.6$.

\includegraphics[bb=10 10 100 120]{sgm42.eps}

\vskip 2mm

\c{${\rm a}=\frac{\mu (u)\rho_M(u)}{2}  <   \pi$ \hskip 40mm ${\rm
a}=\frac{\mu (u)\rho_M(u)}{2}  >  \pi$}\vskip 2mm

\c{\bf Fig.$3.6$}

\vskip 2mm

In an Euclid geometry, a right angle is an angle with measure
$\frac{\pi}{2}$, half of a straight angle and parallel lines are
straight lines never intersecting. They are very important
research objects. Many theorems characterize properties of them in
classical Euclid geometry.

In a map geometry, we can also define a straight angle, a right
angle and parallel $m$-lines by Definition $3.2.2$. Now a {\it
straight angle} is an angle with measure $\pi$ for points not
being vertices of $M$ and $\frac{\rho_M(u)\mu(u)}{2}$ for $\forall
u\in V(M)$. A {\it right angle} is an angle with a half measure of
a straight angle. Two $m$-lines are said {\it parallel} if they
are never intersecting. The following result asserts that map
geometries without boundary are paradoxist geometries.

\vskip 4mm

\no{\bf Theorem $3.2.3$} \ {\it For a map $M$ on a locally
orientable surface with $|M|\geq 3$ and $\rho_M(u)\geq 3$ for
$\forall u\in V(M)$, there exists an angle factor $\mu$ such that
$(M,\mu)$ is a Smarandache geometry by denial the axiom (A5) with
these axioms (A5),(L5) and (R5). }

\vskip 3mm

{\it Proof} \ According to Theorem $3.2.1$, we know that there
exists an angle factor $\mu$ such that there are elliptic
vertices, euclidean vertices and hyperbolic vertices in $(M,\mu)$
simultaneously. The proof is divided into three cases according to
$M$ is planar, orientable or non-orientable. Not loss of
generality, we assume that an angle is measured along a clockwise
direction, i.e., as these cases in Fig.$3.6$ for an $m$-line
passing through an elliptic point or a hyperbolic point.

\vskip 3mm

\no{\bf Case $1.$ \ $M$ is a planar map}

\vskip 2mm

Notice that for a given line $L$ not intersection with the map $M$
and a point $u$ in $(M,\mu)$, if $u$ is an euclidean point, then
there is one and only one line passing through $u$ not
intersecting with $L$, and  if $u$ is an elliptic point, then
there are infinite lines passing through $u$ not intersecting with
$L$, but if $u$ is a hyperbolic point, then each line passing
through $u$ will intersect with $L$. See also in Fig.$3.7$, where
the planar graph is a complete graph $K_4$ and points $1,2$ are
elliptic, the point 3 is euclidean but the point $4$ is
hyperbolic. Then all $m$-lines in the field $A$ do not intersect
with $L$ and each $m$-line passing through the point $4$ will
intersect with the line $L$. Therefore, $(M,\mu)$ is a Smarandache
geometry by denial the axiom (A5) with these axioms (A5), (L5) and
(R5).

\includegraphics[bb=10 10 100 150]{sgm43.eps}

\vskip 2mm

\c{\bf Fig.$3.7$}

\vskip 3mm

\no{\bf Case $2.$ \ $M$ is an orientable map}

\vskip 2mm

According to the classification theorem of compact surfaces, We
only need to prove this result for a torus. Notice that $m$-lines
on a torus has the following property (see [$82$] for details):

\vskip 2mm

{\it If the slope $\varsigma$ of an $m$-line $L$ is a rational
number, then $L$ is a closed line on the torus. Otherwise, $L$ is
infinite, and moreover $L$ passes arbitrarily close to every point
on the torus.}

\vskip 2mm

Whence, if $L_1$ is an $m$-line on a torus with an irrational
slope not passing through an elliptic or a hyperbolic point, then
for any point $u$ exterior to $L_1$, if $u$ is an euclidean point,
then there is only one $m$-line passing through $u$ not
intersecting with $L_1$, and if $u$ is elliptic or hyperbolic, any
$m$-line passing through $u$ will intersect with $L_1$.

Now let $L_2$ be an $m$-line on the torus with an rational slope
not passing through an elliptic or a hyperbolic point, such as the
$m$-line $L_2$ in Fig.$3.8$, $v$ is an euclidean point. If $u$ is
an euclidean point, then each $m$-line $L$ passing through $u$
with rational slope in the area $A$ will not intersect with $L_2$
but each $m$-line passing through $u$ with irrational slope in the
area $A$ will intersect with $L_2$.

\includegraphics[bb=10 10 100 140]{sgm44.eps}

\vskip 3mm

\c{\bf Fig.$3.8$}

\no Therefore, $(M,\mu)$ is a Smarandache geometry by denial the
axiom (A5) with axioms (A5),(L5) and (R5) in this case.

\vskip 3mm

\no{\bf Case $3.$ \ $M$ is a non-orientable map}

\vskip 2mm

Similar to the Case $2$, we only need to prove this result for the
projective plane. An $m$-line in a projective plane is shown in
Fig.$3.9$(a), (b) or (c), in where case (a) is an $m$-line passing
through an euclidean point, (b) passing through an elliptic point
and (c) passing through an hyperbolic point.

\includegraphics[bb=10 10 100 130]{sgm45.eps}

\vskip 2mm

\c{\bf Fig.$3.9$}\vskip 2mm

Now let $L$ be an $m$-line passing through the center in the
circle. Then if $u$ is an euclidean point, there is only one
$m$-line passing through $u$ such as the case $(a)$ in Fig.$3.10$.
If $v$ is an elliptic point then there is an $m$-line passing
through it and intersecting with $L$ such as the case $(b)$ in
Fig.$3.10$. We assume the point $1$ is a point such that there
exists an $m$-line passing through $1$ and $0$, then any $m$-line
in the shade of Fig.$3.10(b)$ passing through $v$ will intersect
with $L$.

\includegraphics[bb=10 10 100 130]{sgm46.eps}

\vskip 2mm

\c{\bf Fig.$3.10$}\vskip 2mm

If $w$ is an euclidean point and there is an $m$-line passing
through it not intersecting with $L$ such as the case $(c)$ in
Fig.$3.10$, then any $m$-line in the shade of Fig.$3.10(c)$
passing through $w$ will not intersect with $L$. Since the
position of the vertices of a map $M$ on a projective plane can be
choose as our wish, we know $(M,\mu)$ is a Smarandache geometry by
denial the axiom (A5) with axioms (A5),(L5) and (R5).

Combining these discussions of Cases $1,2$ and $3$, the proof is
complete. \quad\quad $\natural$.

Similar to Iseri's $s$-manifolds, among map geometries without
boundary there are non-geometries, anti-geometries and
counter-projective geometries, $\cdots$, etc..

\vskip 4mm

\no{\bf Theorem $3.2.4$} \ {\it There are non-geometries in map
geometries without boundary. }

\vskip 3mm

{\it Proof} \ We prove there are map geometries without boundary
satisfying axioms $(A_1^-)-(A_5^-)$. Let $(M,\mu)$ be such a map
geometry with elliptic or hyperbolic points.

($i$) \ Assume $u$ is an eulicdean point and $v$ is an elliptic or
hyperbolic point on $(M,\mu)$. Let $L$ be an $m$-line passing
through points $u$ and $v$ in an Euclid plane. Choose a point $w$
in $L$ after but nearly enough to $v$ when we travel on $L$ from
$u$ to $v$. Then there does not exist a line from $u$ to $w$ in
the map geometry $(M,\mu)$ since $v$ is an elliptic or hyperbolic
point. So the axiom $(A_1^-)$ is true in $(M,\mu)$.

($ii$) \ In a map geometry $(M,\mu)$, an $m$-line maybe closed
such as we have illustrated in the proof of Theorem $3.2.3$.
Choose any two points $A,B$ on a closed $m$-line $L$ in a map
geometry. Then the $m$-line between $A$ and $B$ can not
continuously extend to indefinite in $(M,\mu)$. Whence the axiom
$(A_2^-)$ is true in $(M,\mu)$.

($iii$) \ An $m$-circle in a map geometry is defined to be a set
of continuous points in which all points have a given distance to
a given point. Let $C$ be a $m$-circle in an Euclid plane. Choose
an elliptic or a hyperbolic point $A$ on $C$ which enables us to
get a map geometry $(M,\mu)$. Then $C$ has a gap in $A$ by
definition of an elliptic or hyperbolic point. So the axiom
$(A_3^-)$ is true in a map geometry without boundary.

($iv$) \ By the definition of a right angle, we know that a right
angle on an elliptic point can not equal to a right angle on a
hyperbolic point. So the axiom $(A_4^-)$ is held in a map geometry
with elliptic or hyperbolic points.

($v$) \ The axiom $(A_5^-)$ is true by Theorem $3.2.3$.

Combining these discussions of ($i$)-($v$), we know that there are
non-geometries in map geometries. This completes the proof.
\quad\quad $\natural$ \vskip 3mm

The {\it Hilbert's axiom system} for an Euclid plane geometry
consists five group axioms stated in the following, where we
denote each group by a capital {\it Roman} numeral.

\vskip 3mm

\no{\bf $I.$ Incidence}

\vskip 3mm

\no{\bf $I-1.$} \ {\it For every two points $A$ and $B$, there
exists a line $L$ that contains each of the points $A$ and $B$.}

\no{\bf $I-2.$} \ {\it For every two points $A$ and $B$, there
exists no more than one line that contains each of the points $A$
and $B$.}

\no{\bf $I-3.$} \ {\it There are at least two points on a line.
There are at least three points not on a line.}

\vskip 3mm

\no{\bf $II.$ Betweenness}

\vskip 3mm

\no{\bf $II-1.$} \ {\it If a point $B$ lies between points $A$ and
$C$, then the points $A,B$ and $C$ are distinct points of a line,
and $B$ also lies between $C$ and $A$.}

\no{\bf $II-2.$} \ {\it For two points $A$ and $C$, there always
exists at least one point $B$ on the line $AC$ such that $C$ lies
between $A$ and $B$.}

\no{\bf $II-3.$} \ {\it Of any three points on a line, there
exists no more than one that lies between the other two.}

\no{\bf $II-4.$} \ {\it Let $A,B$ and $C$ be three points that do
not lie on a line, and let $L$ be a line which does not meet any
of the points $A,B$ and $C$. If the line $L$ passes through a
point of the segment $AB$, it also passes through a point of the
segment $AC$, or through a point of the segment $BC$.}

\vskip 3mm

\no{\bf $III.$ Congruence}

\vskip 3mm

\no{\bf $III-1.$} \ {\it If $A_1$ and $B_1$ are two points on a
line $L_1$, and $A_2$ is a point on a line $L_2$ then it is always
possible to find a point $B_2$ on a given side of the line $L_2$
through $A_2$ such that the segment $A_1B_1$ is congruent to the
segment $A_2B_2$.}

\no{\bf $III-2.$} \ {\it If a segment $A_1B_1$ and a segment
$A_2B_2$ are congruent to the segment $AB$, then the segment
$A_1B_1$ is also congruent to the segment $A_2B_2$.}

\no{\bf $III-3.$} \ {\it On the line $L$, let $AB$ and $BC$ be two
segments which except for $B$ have no point in common.
Furthermore, on the same or on another line $L_1$, let $A_1B_1$
and $B_1C_1$ be two segments, which except for $B_1$ also have no
point in common. In that case, if $AB$ is congruent to $A_1B_1$
and $BC$ is congruent to $B_1C_1$, then $AC$ is congruent to
$A_1C_1$.}

\no{\bf $III-4.$} \ {\it Every angle can be copied on a given side
of a given ray in a uniquely determined way.}

\no{\bf $III-5$} \ {\it If for two triangles $ABC$ and
$A_1B_1C_1$, $AB$ is congruent to $A_1B_1$, $AC$ is congruent to
$A_1C_1$ and $\angle BAC$ is congruent to $\angle B_1A_1C_1$, then
$\angle ABC$ is congruent to $\angle A_1B_1C_1$.}

\vskip 3mm

\no{\bf $IV.$ Parallels}

\vskip 3mm

\no{\bf $IV-1.$} \ {\it There is at most one line passes through a
point $P$ exterior a line $L$ that is parallel to $L$.}

\vskip 3mm

\no{\bf $V.$ Continuity}

\vskip 3mm

\no{\bf $V-1$}(Archimedes) \ {\it Let $AB$ and $CD$ be two line
segments with $|AB|\geq |CD|$. Then there is an integer $m$ such
that}

$$m|CD|\leq |AB|\leq (m+1)|CD|.$$

\no{\bf $V-2$}(Cantor) \ {Let $A_1B_1, A_2B_2,\cdots,A_nB_n,
\cdots$ be a segment sequence on a line $L$. If}

$$A_1B_1\supseteq A_2B_2\supseteq\cdots\supseteq A_nB_n\supseteq\cdots,$$

\no{\it then there exists a common point $X$ on each line segment
$A_nB_n$ for any integer $n, n\geq 1$.}

Smarandache defined an anti-geometries by denial some axioms of
Hilbert axiom system for an Euclid geometry. Similar to the
discussion in the reference $[35]$, We obtain the following result
for anti-geometries in map geometries without boundary.

\vskip 4mm

\no{\bf Theorem $3.2.5$} \ {\it Unless axioms $I-3$, $II-3$,
$III-2$, $V-1$ and $V-2$, an anti-geometry can be gotten from map
geometries without boundary by denial other axioms in Hilbert
axiom system.}

\vskip 3mm

{\it Proof} \ The axiom $I-1$ has been denied in the proof of
Theorem $3.2.4$. Since there maybe exists more than one line
passing through two points $A$ and $B$ in a map geometry with
elliptic or hyperbolic points $u$ such as those shown in
Fig.$3.11$. So the axiom $II-2$ can be Smarandachely denied.

\includegraphics[bb=10 10 100 120]{sgm47.eps}

\vskip 2mm

\c{\bf Fig.$3.11$}\vskip 2mm

Notice that an $m$-line maybe has self-intersection points in a
map geometry without boundary. So the axiom $II-1$ can be denied.
By the proof of Theorem $3.2.4$, we know that for two points $A$
and $B$, an $m$-line passing through $A$ and $B$ may not exist.
Whence, the axiom $II-2$ can be denied. For the axiom $II-4$, see
Fig.$3.12$, in where $v$ is a non-euclidean point such that
$\rho_M(v)\mu(v)\geq 2(\pi +\angle ACB)$ in a map geometry.

\includegraphics[bb=10 10 100 120]{sgm48.eps}

\vskip 2mm

\c{\bf Fig.$3.12$}\vskip 2mm

\no So $II-4$ can be also denied. Notice that an $m$-line maybe
has self-intersection points. There are maybe more than one
$m$-lines passing through two given points $A,B$. Therefore, the
axioms $III-1$ and $III-3$ are deniable. For denial the axiom
$III-4$, since an elliptic point $u$ can be measured at most by a
number $\frac{\rho_M(u)\mu(u)}{2} <\pi$, i.e., there is a
limitation for an elliptic point $u$. Whence, an angle with
measure bigger than $\frac{\rho_M(u)\mu(u)}{2}$ can not be copied
on an elliptic point on a given ray.

Because there are maybe more than one $m$-lines passing through
two given points $A$ and $B$ in a map geometry without boundary,
the axiom $III-5$ can be Smarandachely denied in general such as
those shown in Fig.$3.13(a)$ and $(b)$ where $u$ is an elliptic
point.

\includegraphics[bb=10 10 100 120]{sgm49.eps}

\vskip 2mm

\c{\bf Fig.$3.13$}\vskip 2mm

For the parallel axiom $IV-1$, it has been denied by the proofs of
Theorems $3.2.3$ and $3.2.4$.

Notice that axioms $I-3$, $II-3$ $III-2$, $V-1$ and $V-2$ can not
be denied in a map geometry without boundary. This completes the
proof. \quad\quad $\natural$

For counter-projective geometries, we have a result as in the
following.

\vskip 4mm

\no{\bf Theorem $3.2.6$} \ {\it Unless the axiom $(C3)$, a
counter-projective geometry can be gotten from map geometries
without boundary by denial axioms $(C1)$ and $(C2)$.}

\vskip 3mm

{\it Proof} \  Notice that axioms $(C1)$ and $(C2)$ have been
denied in the proof of Theorem $3.2.5$. Since a map is embedded on
a locally orientable surface, every $m$-line in a map geometry
without boundary may contains infinite points. Therefore the axiom
$(C3)$ can not be Smarandachely denied. \quad\quad $\natural$

\vskip 5mm

\no{\bf \S $3.3$ \ Map Geometries with Boundary}

\vskip 4mm

\no A {\it Poincar\'{e}'s model} for a hyperbolic geometry is an
upper half-plane in which lines are upper half-circles with center
on the $x$-axis or upper straight lines perpendicular to the
$x$-axis such as those shown in Fig.$3.14$.

\includegraphics[bb=10 10 100 130]{sgm50.eps}

\vskip 2mm

\c{\bf Fig.$3.14$}\vskip 2mm

If we think that all infinite points are the same, then a
Poincar\'{e}'s model for a hyperbolic geometry is turned to a {\it
Klein model} for a hyperbolic geometry which uses a boundary
circle and lines are  straight line segment in this circle, such
as those shown in Fig.$3.15$.

\includegraphics[bb=0 10 100 120]{sgm51.eps}

\vskip 3mm

\c{\bf Fig.$3.15$}\vskip 2mm

By a combinatorial map view, a Klein's model is nothing but a one
face map geometry. This fact hints us to introduce map geometries
with boundary, which is defined in the next definition.

\vskip 4mm

\no{\bf Definition $3.3.1$} \ {\it For a map geometry $(M,\mu )$
without boundary and faces $f_1,f_2,\cdots ,f_l$ $\in F(M), 1\leq
l\leq \phi (M)-1$, if $S(M)\setminus\{f_1,f_2,\cdots ,f_l\}$ is
connected, then call $(M,\mu)^{-l}= (S(M)\setminus\{f_1,f_2,\cdots
,f_l\}, \mu)$ a map geometry with boundary $f_1,f_2,\cdots ,f_l$
and orientable or not if $(M,\mu)$ is orientable or not, where
$S(M)$ denotes the locally orientable surface on which $M$ is
embedded.}

\vskip 3mm

The $m$-points and $m$-lines in a map geometry $(M,\mu)^{-l}$ are
defined as same as Definition $3.2.3$ by adding an $m$-line
terminated at the boundary of this map geometry. Two $m^-$-lines
on the torus and projective plane are shown in these Fig.$3.16$
and Fig.$3.17$, where the shade field denotes the boundary.

\includegraphics[bb=10 10 100 125]{sgm52.eps}

\vskip 2mm

\c{\bf Fig.$3.16$}

\includegraphics[bb=10 10 100 135]{sgm53.eps}

\vskip 2mm

\c{\bf Fig.$3.17$}

\vskip 2mm

All map geometries with boundary are also Smarandache geometries
which is convince by a result in the following.

\vskip 4mm

\no{\bf Theorem $3.3.1$} {\it For a map $M$ on a locally
orientable surface with order$\geq 3$, vertex valency$\geq 3$ and
a face $f\in F(M)$, there is an angle factor $\mu$ such that
$(M,\mu)^{-1}$ is a Smarandache geometry by denial the axiom (A5)
with these axioms (A5),(L5) and (R5). }

\vskip 3mm

{\it Proof} \ Similar to the proof of Theorem $3.2.3$, we consider
a map $M$ being a planar map, an orientable map on a torus or a
non-orientable map on a projective plane, respectively. We can get
the assertion. In fact, by applying the property that $m$-lines in
a map geometry with boundary are terminated at the boundary, we
can get an more simpler proof for this theorem. \quad\quad
$\natural$

Notice that in a one face map geometry $(M,\mu)^{-1}$ with
boundary is just a Klein's model for hyperbolic geometry if we
choose all points being euclidean.

Similar to map geometries without boundary, we can also get
non-geometries, anti-geometries and counter-projective geometries
from map geometries with boundary.

\vskip 4mm

\no{\bf Theorem $3.3.2$} \ {\it There are non-geometries in map
geometries with boundary. }

\vskip 3mm

{\it Proof} \ The proof is similar to the proof of Theorem $3.2.4$
for map geometries without boundary. Each of axioms
$(A_1^-)-(A_5^-)$ is hold, for example, cases $(a)-(e)$ in
Fig.$3.18$,

\includegraphics[bb=10 10 100 230]{sgm54.eps}

\vskip 2mm

\c{\bf Fig.$3.18$}\vskip 3mm

\no in where there are no an $m$-line from points $A$ to $B$ in
$(a)$, the line $AB$ can not be continuously extended to
indefinite in $(b)$, the circle has gap in $(c)$, a right angle at
an euclidean point $v$ is not equal to a right angle at an
elliptic point $u$ in $(d)$ and there are infinite $m$-lines
passing through a point $P$ not intersecting with the $m$-line $L$
in $(e)$. Whence, there are non-geometries in map geometries with
boundary. \quad\quad $\natural$

\vskip 4mm

\no{\bf Theorem $3.3.3$} \ {\it Unless axioms $I-3$, $II-3$
$III-2$, $V-1$ and $V-2$ in the Hilbert's axiom system for an
Euclid geometry, an anti-geometry can be gotten from map
geometries with boundary by denial other axioms in this axiom
system.}

\vskip 4mm

\no{\bf Theorem $3.3.4$} \ {\it Unless the axiom $(C3)$, a
counter-projective geometry can be gotten from map geometries with
boundary by denial axioms $(C1)$ and $(C2)$.}

\vskip 3mm

{\it Proof} \ The proofs of Theorems $3.3.3$ and $3.3.4$ are
similar to the proofs of Theorems $3.2.5$ and $3.2.6$. The reader
is required to complete their proof. \quad\quad $\natural$

\vskip 5mm

\no{\bf $\S 3.4$ \ The Enumeration of Map Geometries}

\vskip 4mm

\no For classifying map geometries, the following definition is
needed.

\vskip 4mm

\no{\bf Definition $3.4.1$} \ {\it Two map geometries
$(M_1,\mu_1)$ and $(M_2,\mu_2)$ or $(M_1,\mu_1)^{-l}$ and
$(M_2,\mu_2)^{-l}$ are said to be equivalent each other if there
is a bijection $\theta :M_1\rightarrow M_2$ such that for $\forall
u\in V(M)$, $\theta (u)$ is euclidean, elliptic or hyperbolic if
and only if $u$ is euclidean, elliptic or hyperbolic.}

\vskip 3mm

A relation for the numbers of unrooted maps with map geometries is
in the following result.

\vskip 4mm

\no{\bf Theorem $3.4.1$} \ {\it Let ${\mathcal M}$ be a set of
non-isomorphic maps of order $n$ and with $m$ faces. Then the
number of map geometries without boundary is $3^n|\mathcal M|$ and
the number of map geometries with one face being its boundary is
$3^nm|\mathcal M|$.}

\vskip 3mm

{\it Proof} \ By the definition of equivalent map geometries, for
a given map $M\in {\mathcal M}$, there are $3^n$ map geometries
without boundary and $3^nm$ map geometries with one face being its
boundary by Theorem $3.3.1$. Whence, we get $3^n|\mathcal M|$ map
geometries without boundary and $3^nm|\mathcal M|$ map geometries
with one face being its boundary from ${\mathcal M}$.\quad\quad
$\natural$.

We get an enumeration result for non-equivalent map geometries
without boundary as follows.

\vskip 4mm

\no{\bf Theorem $3.4.2$} \ {\it The numbers $n^O(\Gamma ,g)$ and
$n^N(\Gamma ,g)$ of non-equivalent orientable and non-orientable
map geometries without boundary underlying a simple graph $\Gamma$
by denial the axiom (A5) by (A5), (L5) or (R5) are  }

$$n^O(\Gamma ,g)=\frac{3^{|\Gamma|}\prod\limits_{v\in V(\Gamma)}
(\rho (v)-1)!}{2|{\rm Aut}\Gamma |},$$

\no{\it and}

$$n^N(\Gamma ,g)=\frac{(2^{\beta (\Gamma )}-1)3^{|\Gamma|}
\prod\limits_{v\in V(\Gamma)}(\rho (v)-1)!}{2|{\rm Aut}\Gamma
|},$$

\no{\it where $\beta (\Gamma)=\varepsilon (\Gamma)-\nu (\Gamma)
+1$ is the Betti number of the graph $\Gamma$.}

\vskip 3mm

{\it Proof} \ Denote the set of non-isomorphic maps underlying the
graph $\Gamma$ on locally orientable surfaces by ${\mathcal
M}(\Gamma)$ and the set of embeddings of the graph $\Gamma$ on
locally orientable surfaces by ${\mathcal E}(\Gamma)$. For a map
$M, M\in {\mathcal M}(\Gamma)$, there are $\frac{3^{|M|}}{|{\rm
Aut}M|}$ different map geometries without boundary by choice the
angle factor $\mu$ on a vertex $u$ such that $u$ is euclidean,
elliptic or hyperbolic. From permutation groups, we know that

$$
|{\rm Aut}\Gamma\times \left<\alpha\right>|=|({\rm
Aut}\Gamma)_M||M^{{\rm Aut}\Gamma\times
\left<\alpha\right>}|=|{\rm Aut}M||M^{{\rm Aut}\Gamma\times
\left<\alpha\right>}|.
$$

\no Therefore, we get that

\begin{eqnarray*}
n^O(\Gamma ,g)& = & \sum\limits_{M\in {\mathcal
M}(\Gamma)}\frac{3^{|M|}}{|{\rm Aut}M|}\\
&=& \frac{3^{|\Gamma |}}{|\rm Aut\Gamma\times
\left<\alpha\right>|}\sum\limits_{M\in {\mathcal
M}(\Gamma)}\frac{|\rm
Aut\Gamma\times \left<\alpha\right>|}{|{\rm Aut}M|}\\
&=& \frac{3^{|\Gamma |}}{|\rm Aut\Gamma\times
\left<\alpha\right>|}\sum\limits_{M\in {\mathcal
M}(\Gamma)}|M^{{\rm
Aut}\Gamma\times \left<\alpha\right>}|\\
&=& \frac{3^{|\Gamma |}}{|\rm Aut\Gamma\times
\left<\alpha\right>|}|{\mathcal
E}^O(\Gamma)|\\
&=& \frac{ 3^{|\Gamma|}\prod\limits_{v\in V(\Gamma)}(\rho
(v)-1)!}{2|{\rm Aut}\Gamma |}.
\end{eqnarray*}

\no Similarly, we can also get that

\begin{eqnarray*}
n^N(\Gamma ,g)& = & \frac{3^{|\Gamma |}}{|\rm Aut\Gamma\times
\left<\alpha\right>|}|{\mathcal
E}^N(\Gamma)|\\
&=& \frac{ (2^{\beta (\Gamma )}-1)3^{|\Gamma|}\prod\limits_{v\in
V(\Gamma)}(\rho (v)-1)!}{2|{\rm Aut}\Gamma |}.
\end{eqnarray*}

\no This completes the proof. \quad\quad $\natural$

For classifying map geometries with boundary, we get a result as
in the following.

\vskip 4mm

\no{\bf Theorem $3.4.3$} \ {\it The numbers $n^O(\Gamma ,-g)$,
$n^N(\Gamma ,-g)$ of non-equivalent orientable, non-orientable map
geometries with one face being its boundary underlying a simple
graph $\Gamma$ by denial the axiom (A5) by (A5), (L5) or (R5) are
respective}

$$n^O(\Gamma ,-g)=\frac{3^{|\Gamma|}}{2|{\rm Aut}\Gamma|}[(\beta (\Gamma)+1)
\prod\limits_{v\in V(\Gamma)}(\rho (v)-1)!-\frac{2d(g[\Gamma
](x))}{dx}|_{x=1}]$$

\no{\it and}

$$n^N(\Gamma ,-g)=\frac{(2^{\beta (\Gamma)}-1)3^{|\Gamma|}}{2|{\rm Aut}\Gamma|}[(\beta (\Gamma)+1)
\prod\limits_{v\in V(\Gamma)}(\rho (v)-1)!-\frac{2d(g[\Gamma
](x))}{dx}|_{x=1}],$$

\no{\it where $g[\Gamma ](x)$ is the genus polynomial of the graph
$\Gamma$, i.e., $g[\Gamma ](x)=\sum\limits_{k=\gamma
(\Gamma)}^{\gamma_m(\Gamma)}g_k[\Gamma ]x^k$ with $g_k[\Gamma ]$
being the number of embeddings of $\Gamma$ on the orientable
surface of genus $k$.}

\vskip 3mm

{\it Proof} \ Notice that $\nu (M)-\varepsilon (M)+\phi
(M)=2-2g(M)$ for an orientable map $M$ by the
Euler-Poincar$\acute{e}$ formula. Similar to the proof of Theorem
$3.4.2$ with the same meaning for ${\mathcal M}(\Gamma)$, we know
that

\begin{eqnarray*}
n^O(\Gamma ,-g)&=& \sum\limits_{M\in{\mathcal
M}(\Gamma)}\frac{\phi (M)3^{|M|}}{|{\rm Aut}M|}\\
&=&\sum\limits_{M\in{\mathcal M}(\Gamma)}\frac{(2+\varepsilon
(\Gamma)-\nu (\Gamma)-2g(M))3^{|M|}}{|{\rm Aut}M|}\\
&=&\sum\limits_{M\in{\mathcal M}(\Gamma)}\frac{(2+\varepsilon
(\Gamma)-\nu (\Gamma))3^{|M|}}{|{\rm
Aut}M|}-\sum\limits_{M\in{\mathcal
M}(\Gamma)}\frac{2g(M)3^{|M|}}{|{\rm Aut}M|}\\
&=&\frac{(2+\varepsilon (\Gamma)-\nu (\Gamma))3^{|M|}}{|{\rm
Aut}\Gamma\times \left<\alpha\right>|}\sum\limits_{M\in{\mathcal
M}(\Gamma)}\frac{|{\rm Aut}\Gamma\times
\left<\alpha\right>|}{|{\rm
Aut}M|}\\
&-&\frac{2\times 3^{|\Gamma |}}{|{\rm Aut}\Gamma\times
\left<\alpha\right>|}\sum\limits_{M\in{\mathcal
M}(\Gamma)}\frac{g(M)|{\rm Aut}\Gamma\times
\left<\alpha\right>|}{|{\rm Aut}M|}\\
&=& \frac{(\beta (\Gamma)+1)3^{|M|}}{|{\rm Aut}\Gamma\times
\left<\alpha\right>|}\sum\limits_{M\in{\mathcal
M}}(\Gamma)|M^{{\rm
Aut}\Gamma\times \left<\alpha\right>}|\\
&-&\frac{3^{|\Gamma |}}{|{\rm Aut}\Gamma
|}\sum\limits_{M\in{\mathcal M}(\Gamma)}g(M)|M^{{\rm
Aut}\Gamma\times \left<\alpha\right>}|\\
&=&\frac{(\beta (\Gamma)+1)3^{|\Gamma |}}{2|{\rm Aut}\Gamma
|}\prod\limits_{v\in V(\Gamma)}(\rho (v)-1)! -\frac{3^{|\Gamma
|}}{|{\rm Aut}\Gamma |}\sum\limits_{k=\gamma
(\Gamma)}^{\gamma_m(\Gamma)}kg_k[\Gamma ]\\
&=&\frac{3^{|\Gamma |}}{2|{\rm Aut}\Gamma |}[(\beta
(\Gamma)+1)\prod\limits_{v\in V(\Gamma)}(\rho (v)-1)!-
\frac{2d(g[\Gamma ](x))}{dx}|_{x=1}].
\end{eqnarray*}

\no by Theorem $3.4.1$.

Notice that $n^L(\Gamma ,-g)=n^O(\Gamma ,-g)+n^N(\Gamma ,-g)$ and
the number of re-embeddings an orientable map $M$ on surfaces is
$2^{\beta (M)}$ (see also [$56$] for details). We know that

\begin{eqnarray*}
n^L(\Gamma ,-g)&=&\sum\limits_{M\in {\mathcal
M}(\Gamma)}\frac{2^{\beta (M)}\times 3^{|M|}\phi (M)}{|{\rm
Aut}M|}\\
&=& 2^{\beta (M)}n^O(\Gamma ,-g).
\end{eqnarray*}

\no Whence, we get that

\begin{eqnarray*}
n^N(\Gamma ,-g)&=& (2^{\beta (M)}-1)n^O(\Gamma ,-g)\\
&=& \frac{(2^{\beta (M)}-1)3^{|\Gamma |}}{2|{\rm Aut}\Gamma
|}[(\beta (\Gamma)+1)\prod\limits_{v\in V(\Gamma)}(\rho (v)-1)!-
\frac{2d(g[\Gamma ](x))}{dx}|_{x=1}].
\end{eqnarray*}

\no This completes the proof.\quad\quad $\natural$

\vskip 5mm

\no{\bf \S $3.5$ \ Remarks and Open Problems}

\vskip 4mm

\no{\bf $3.5.1$} \ A complete Hilbert axiom system for an Euclid
geometry contains axioms $I-i, 1\leq i\leq 8$; $II-j, 1\leq j\leq
4$; $III-k, 1\leq k\leq 5$; $IV-1$ and $V-l, 1\leq l\leq 2$, which
can be also applied to the geometry of space. Unless $I-i, 4\leq
i\leq 8$, other axioms are presented in Section $3.2$. Each of
axioms $I-i,4\leq i\leq 8$ is described in the following.

\vskip 3mm

\no{\bf $I-4$} \ {\it For three non-collinear points $A,B$ and
$C$, there is one and only one plane passing through them.}

\no{\bf $I-5$} \ {\it Each plane has at least one point.}

\no{\bf $I-6$} \ {\it If two points $A$ and $B$ of a line $L$ are
in a plane $\sum$, then every point of $L$ is in the plane
$\sum$.}

\no{\bf $I-7$} \ {\it If two planes $\sum_1$ and $\sum_2$ have a
common point $A$, then they have another common point $B$.}

\no{\bf $I-8$} \ {\it There are at least four points not in one
plane.}\vskip 2mm

By the Hilbert's axiom system, the following result for parallel
planes can be obtained.

\vskip 3mm

({\bf T}) \ {\it Passing through a given point $A$ exterior to a
given plane $\sum$ there is one and only one plane parallel to
$\sum$.}

\vskip 2mm

This result seems like the Euclid's fifth axiom. Similar to the
Smarandache's notion, we present problems by denial this theorem
for the geometry of space as follows.

\vskip 4mm

\no{\bf Problem $3.5.1$} \ {\it Construct a geometry of space by
denial the parallel theorem of planes with}

\vskip 3mm

({\bf $T_1^-$}) \ {\it there are at least a plane $\sum$ and a
point $A$ exterior to the plane $\sum$ such that no parallel plane
to $\sum$ passing through the point $A$.}

({\bf $T_2^-$}) \ {\it there are at least a plane $\sum$ and a
point $A$ exterior to the plane $\sum$ such that there are finite
parallel planes to $\sum$ passing through the point $A$.}

({\bf $T_3^-$}) \ {\it there are at least a plane $\sum$ and a
point $A$ exterior to the plane $\sum$ such that there are
infinite parallel planes to $\sum$ passing through the point $A$.}

\vskip 4mm

\no{\bf Problem $3.5.2$} \ {\it Similar to the Iseri's idea define
an elliptic, euclidean, or hyperbolic point or plane in ${\bf
R}^3$ and apply these Plato polyhedrons to construct Smarandache
geometries of a space ${\bf R}^3$.}

\vskip 4mm

\no{\bf Problem $3.5.3$} \ {\it Similar to map geometries define
graph in a space geometries and apply graphs in ${\bf R}^3$ to
construct Smarandache geometries of a space ${\bf R}^3$.}

\vskip 4mm

\no{\bf Problem $3.5.4$} \ {\it For an integer $n, n\geq 4$,
define Smarandache geometries in ${\bf R}^n$ by denial some axioms
for an Euclid geometry in ${\bf R}^n$ and construct them.}

\vskip 3mm

\no{\bf $3.5.2$} \ The terminology {\it map geometry} was first
appeared in $[55]$ which enables us to find non-homogenous spaces
from already known homogenous spaces and is also a typical example
for application combinatorial maps to metric geometries. Among
them there are many problems not solved yet until today. Here we
would like to describe some of them.

\vskip 4mm

\no{\bf Problem $3.5.5$} \ {\it For a given graph $G$, determine
non-equivalent map geometries with an underlying graph $G$,
particularly, for graphs $K_n$, $K(m,n), m,n\geq 4$ and enumerate
them.}

\vskip 4mm

\no{\bf Problem $3.5.6$} \ {\it For a given locally orientable
surface $S$, determine non-equivalent map geometries on $S$, such
as a sphere, a torus or a projective plane, $\cdots$ and enumerate
them.}

\vskip 4mm

\no{\bf Problem $3.5.7$} \ {\it Find characteristics for
equivalent map geometries or establish new ways for classifying
map geometries.}

\vskip 4mm

\no{\bf Problem $3.5.8$} \ {\it Whether can we rebuilt an
intrinsic geometry on surfaces, such as a sphere, a torus or a
projective plane, $\cdots$, by map geometries?}

\newpage

\no{\large\bf $4.$ Planar map geometries}

\vskip 15mm

\no Fundamental elements in an Euclid geometry are those of
points, lines, polygons and circles. For a map geometry, the
situation is more complex since a point maybe an elliptic,
euclidean or a hyperbolic point, a polygon maybe a line, $\cdots$,
etc.. This chapter concentrates on discussing fundamental elements
and measures such as angle, area, curvature, $\cdots$, etc., also
parallel bundles in planar map geometries, which can be seen as a
first step for comprehending map geometries on surfaces. All
materials of this chapter will be used in Chapters $5$-$6$ for
establishing relations of an integral curve with a differential
equation system in a pseudo-plane geometry and continuous
phenomena with discrete phenomena

\vskip 5mm

\no{\bf $\S 4.1$ \ Points in a Planar Map Geometry}

\vskip 4mm

\no Points in a map geometry are classified into three classes:
{\it elliptic, euclidean} and { \it hyperbolic}. There are only
finite non-euclidean points considered in Chapter $3$ because we
had only defined an elliptic, euclidean or a hyperbolic point on
vertices of a map. In a planar map geometry, we can present an
even more delicate consideration for euclidean or non-euclidean
points and find infinite non-euclidean points in a plane.

Let $(M,\mu)$ be a planar map geometry on a plane $\sum$. Choose
vertices $u,v\in V(M)$. A mapping is called an {\it angle function
between $u$ and $v$} if there is a smooth monotone mapping $f: \
\sum\rightarrow \sum$ such that $f(u)=\frac{\rho_M(u)\mu(u)}{2}$
and $f(v)=\frac{\rho_M(v)\mu(v)}{2}$. Not loss of generality, we
can assume that each edge in a planar map geometry is an angle
function. Then we know a result as in the following.

\vskip 4mm

\no{\bf Theorem $4.1.1$} \ {\it A planar map geometry $(M,\mu)$
has infinite non-euclidean points if and only if there is an edge
$e=(u,v)\in E(M)$ such that $\rho_M(u)\mu(u)\not=\rho_M(v)\mu(v)$,
or $\rho_M(u)\mu(u)$ is a constant but $\not= 2\pi$ for $\forall
u\in V(M)$, or a loop $(u,u)\in E(M)$ attaching a non-euclidean
point $u$.}

\vskip 3mm

{\it Proof} \ If there is an edge $e=(u,v)\in E(M)$ such that
$\rho_M(u)\mu(u)\not=\rho_M(v)\mu(v)$, then at least one of
vertices $u$ and $v$ in $(M,\mu)$ is non-euclidean. Not loss of
generality, we assume the vertex $u$ is non-euclidean.

If $u$ and $v$ are elliptic or $u$ is elliptic but $v$ is
euclidean, then by the definition of angle functions, the edge
$(u,v)$ is correspondent with an angle function $f:
\sum\rightarrow \sum$ such that $f(u)=\frac{\rho_M(u)\mu(u)}{2}$
and $f(v)=\frac{\rho_M(v)\mu(v)}{2}$, each points is non-euclidean
in $(u,v)\setminus\{v\}$. If $u$ is elliptic but $v$ is
hyperbolic, i.e., $\rho_M(u)\mu(u) < 2\pi$ and $\rho_M(v)\mu(v) >
2\pi$, since $f$ is smooth and monotone on $(u,v)$, there is one
and only one point $x^*$ in $(u,v)$ such that $f(x^*)=\pi$.
Thereby there are infinite non-euclidean points on $(u,v)$.

Similar discussion can be gotten for the cases that $u$ and $v$
are both hyperbolic, or $u$ is hyperbolic but $v$ is euclidean, or
$u$ is hyperbolic but $v$ is elliptic.

If $\rho_M(u)\mu(u)$ is a constant but $\not=2\pi$ for $\forall
u\in V(M)$, then each point on an edges is a non-euclidean point.
Thereby there are infinite non-euclidean points in $(M,\mu)$.

Now if there is a loop $(u,u)\in E(M)$ and $u$ is non-eucliean,
then by definition, each point $v$ on the loop $(u,u)$ satisfying
that $f(v)>$ or $< \pi$ according to $\rho_M(u)\mu(u)> \pi$ or $<
\pi$. Therefore there are also infinite non-euclidean points on
the loop $(u,u)$.

On the other hand, if there are no an edge $e=(u,v)\in E(M)$ such
that $\rho_M(u)\mu(u)\not=\rho_M(v)\mu(v)$, i.e.,
$\rho_M(u)\mu(u)=\rho_M(v)\mu(v)$ for $\forall (u,v)\in E(M)$, or
there are no vertices $ u\in V(M)$ such that $\rho_M(u)\mu(u)$ is
a constant but $\not=2\pi$ for $\forall$, or there are no loops
$(u,u)\in E(M)$ with a non-eucliean point $u$, then all angle
functions on these edges of $M$ are an constant $\pi$. Therefore
there are no non-euclidean points in the map geometry $(M,\mu)$.
This completes the proof.\quad\quad $\natural$

For euclidean points in a planar map geometry $(M,\mu)$, we get
the following result.

\vskip 4mm

\no{\bf Theorem $4.1.2$} \ {\it For a planar map geometry
$(M,\mu)$ on a plane $\sum$,

($i$) \ every point in $\sum\setminus E(M)$ is an euclidean point;

($ii$) \ there are infinite euclidean points on $M$ if and  only
if there exists an edge $(u,v)\in E(M)$ ($u=v$ or $u\not= v$) such
that $u$ and $v$ are both euclidean.}

\vskip 3mm

{\it Proof} \ By the definition of angle functions, we know that
every point is euclidean if it is not on $M$. So the assertion
($i$) is true.

According to the proof of Theorem $4.1.1$, there are only finite
euclidean points unless there is an edge $(u,v)\in E(M)$ with
$\rho_M(u)\mu(u)=\rho_M(v)\mu(v)=2\pi$. In this case, there are
infinite euclidean points on the edge $(u,v)$. Thereby the
assertion ($ii$) is also holds.\quad\quad $\natural$

According to Theorems $4.1.1$ and $4.1.2$, we classify edges in a
planar map geometry $(M,\mu)$ into six classes as follows.

\vskip 3mm

{\bf $C_E^1$ (euclidean-elliptic edges)}: \ {\it edges $(u,v)\in
E(M)$ with $\rho_M(u)\mu(u) = 2\pi$ but $\rho_M(v)\mu(v) < 2\pi$.}

{\bf $C_E^2$ (euclidean-euclidean edges)}: \ {\it edges $(u,v)\in
E(M)$ with $\rho_M(u)\mu(u) = 2\pi$ and $\rho_M(v)\mu(v) = 2\pi$.}

{\bf $C_E^3$ (euclidean-hyperbolic edges)}: \ {\it edges $(u,v)\in
E(M)$ with $\rho_M(u)\mu(u) = 2\pi$ but $\rho_M(v)\mu(v) > 2\pi$.}

{\bf $C_E^4$ (elliptic-elliptic edges)}: \ {\it edges $(u,v)\in
E(M)$ with $\rho_M(u)\mu(u) < 2\pi$ and $\rho_M(v)\mu(v) < 2\pi$.}

{\bf $C_E^5$ (elliptic-hyperbolic edges)}: \ {\it edges $(u,v)\in
E(M)$ with $\rho_M(u)\mu(u) < 2\pi$ but $\rho_M(v)\mu(v) > 2\pi$.}

{\bf $C_E^6$ (hyperbolic-hyperbolic edges)}: \ {\it edges
$(u,v)\in E(M)$ with $\rho_M(u)\mu(u) > 2\pi$ and $\rho_M(v)\mu(v)
> 2\pi$.}

\vskip 2mm

In Fig.$4.1(a)-(f)$, these $m$-lines passing through an edge in
one of classes of $C_E^1$-$C_E^6$ are shown, where $u$ is elliptic
and $v$ is eucildean in $(a)$, $u$ and $v$ are both euclidean in
$(b)$, $u$ is eucildean but $v$ is hyperbolic in $(c)$, $u$ and
$v$ are both elliptic in (d), $u$ is elliptic but $v$ is
hyperbolic in $(e)$ and $u$ and $v$ are both hyperbolic in $(f)$,
respectively.

\includegraphics[bb=10 10 400 220]{sgm55.eps}

\vskip 2mm

\c{\bf Fig.$4.1$}

\vskip 2mm

Denote by $V_{el}(M), V_{eu}(M)$ and $V_{hy}(M)$ the respective
sets of elliptic, euclidean and hyperbolic points in $V(M)$ in a
planar map geometry $(M,\mu)$. Then we get a result as in the
following.

\vskip 4mm

\no{\bf Theorem $4.1.3$} \ {\it Let $(M,\mu)$ be a planar map
geometry. Then}

$$
\sum\limits_{u\in V_{el}(M)}\rho_M(u)+\sum\limits_{v\in
V_{eu}(M)}\rho_M(v)+\sum\limits_{w\in V_{hy}(M)}\rho_M(w)
=2\sum\limits_{i=1}^6|C_E^i|
$$

\no{\it and}

$$|V_{el}(M)|+|V_{eu}(M)|+|V_{hy}(M)|+\phi(M)=\sum\limits_{i=1}^6|C_E^i|+2.$$

\no{\it where $\phi(M)$ denotes the number of faces of a map $M$.}

\vskip 3mm

{\it Proof} \ Notice that

$$|V(M)|=|V_{el}(M)|+|V_{eu}(M)|+|V_{hy}(M)| \ {\rm and} \ |E(M)|=\sum\limits_{i=1}^6|C_E^i|$$

\no for a planar map geometry $(M,\mu)$. By two well-known results

$$\sum\limits_{v\in V(M)}\rho_M(v)= 2|E(M)| \ \ {\rm and} \ \ |V(M)|-|E(M)|+\phi(M)=2$$

\no for a planar map $M$, we know that

$$
\sum\limits_{u\in V_{el}(M)}\rho_M(u)+\sum\limits_{v\in
V_{eu}(M)}\rho_M(v)+\sum\limits_{w\in V_{hy}(M)}\rho_M(w)
=2\sum\limits_{i=1}^6|C_E^i|
$$

\no and

$$|V_{el}(M)|+|V_{eu}(M)|+|V_{hy}(M)|+\phi(M)=\sum\limits_{i=1}^6|C_E^i|+2.\ \ \natural$$

\vskip 5mm

\no{\bf \S $4.2$ \ Lines in a Planar Map Geometry}

\vskip 4mm

\no The situation of $m$-lines in a planar map geometry $(M,\mu)$
is more complex. Here an $m$-line maybe open or closed, with or
without self-intersections in a plane. We discuss all of these
$m$-lines and their behaviors in this section, .

\vskip 4mm

\no{\bf $4.2.1.$ Lines in a planar map geometry}

\vskip 3mm

\no As we have seen in Chapter $3$, $m$-lines in a planar map
geometry $(M,\mu)$ can be classified into three classes.

\vskip 3mm

{\bf $C_L^1$(opened lines without self-intersections)}: \ {\it
$m$-lines in $(M,\mu)$ have an infinite number of continuous
$m$-points without self-intersections and endpoints and may be
extended indefinitely in both directions.}

\vskip 2mm

{\bf $C_L^2$(opened lines with self-intersections)}: \ {\it
$m$-lines in $(M,\mu)$ have an infinite number of continuous
$m$-points and self-intersections but without endpoints and may be
extended indefinitely in both directions.}

\vskip 2mm

{\bf $C_L^3$(closed lines)}: \ {\it $m$-lines in $(M,\mu)$ have an
infinite number of continuous $m$-points and will come back to the
initial point as we travel along any one of these $m$-lines
starting at an initial point.}

\vskip 2mm

By this classification, a straight line in an Euclid plane is
nothing but an opened $m$-line without non-euclidean points.
Certainly, $m$-lines in a planar map geometry $(M,\mu)$ maybe
contain non-euclidean points. In Fig.$4.2$, these $m$-lines shown
in $(a),(b)$ and $(c)$ are opened $m$-line without
self-intersections, opened $m$-line with a self-intersection and
closed $m$-line with $A,B,C,D$ and $E$ non-euclidean points,
respectively.

\includegraphics[bb=10 10 400 110]{sgm56.eps}

\vskip 3mm

\c{\bf Fig.$4.2$}

\vskip 2mm

Notice that a closed $m$-line in a planar map geometry maybe also
has self-intersections. A closed $m$-line is said to be {\it
simply closed} if it has no self-intersections, such as the
$m$-line in Fig.$4.2(c)$. For simply closed $m$-lines, we know the
following result.

\vskip 4mm

\no{\bf Theorem $4.2.1$} \ {\it Let $(M,\mu)$ be a planar map
geometry. An $m$-line $L$ in $(M,\mu)$ passing through $n$
non-euclidean points $x_1,x_2,\cdots, x_n$ is simply closed if and
only if}

$$\sum\limits_{i=1}^nf(x_i)=(n-2)\pi,$$

\no{\it where $f(x_i)$ denotes the angle function value at an
$m$-point $x_i, 1\leq i\leq n$.}

\vskip 3mm

{\it Proof} \ By results in an Euclid geometry of plane, we know
that the angle sum of an $n$-polygon is $(n-2)\pi$. In a planar
map geometry $(M,\mu)$, a simply closed $m$-line $L$ passing
through $n$ non-euclidean points $x_1,x_2,\cdots, x_n$ is nothing
but an $n$-polygon with vertices $x_1,x_2,\cdots, x_n$. Whence, we
get that

$$\sum\limits_{i=1}^nf(x_i)=(n-2)\pi.$$

Now if a simply $m$-line $L$ passing through $n$ non-euclidean
points $x_1,x_2,\cdots, x_n$ with

$$\sum\limits_{i=1}^nf(x_i)=(n-2)\pi$$

\no held, then $L$ is nothing but an $n$-polygon with vertices
$x_1,x_2,\cdots, x_n$. Therefore, $L$ is simply closed.\quad\quad
$\natural$

By applying Theorem $4.2.1$, we can also find conditions for an
opened $m$-line with or without self-intersections.

\vskip 4mm

\no{\bf Theorem $4.2.2$} \ {\it Let $(M,\mu)$ be a planar map
geometry. An $m$-line $L$ in $(M,\mu)$ passing through $n$
non-euclidean points $x_1,x_2,\cdots, x_n$ is opened without
self-intersections if and only if $m$-line segments $x_ix_{i+1},
1\leq i\leq n-1$ are not intersect two by two and}

$$\sum\limits_{i=1}^nf(x_{i})\geq (n-1)\pi.$$

\vskip 3mm

{\it Proof} \ By the Euclid's fifth postulate for a plane
geometry, two straight lines will meet on the side on which the
angles less than two right angles if we extend them to
indefinitely. Now for an $m$-line $L$ in a planar map geometry
$(M,\mu)$, if it is opened without self-intersections, then for
any integer $i, 1\leq i\leq n-1$, $m$-line segments $x_{i}x_{i+1}$
will not intersect two by two and the $m$-line $L$ will also not
intersect before it enters $x_1$ or leaves $x_n$.

\includegraphics[bb=10 10 400 105]{sgm57.eps}

\vskip 2mm

\c{\bf Fig.$4.3$}

\vskip 2mm

Now look at Fig.$4.3$, in where line segment $x_1x_n$ is an added
auxiliary $m$-line segment. We know that

$$\angle 1+\angle 2= f(x_1) \ {\rm and} \ \angle 3+\angle 4=f(x_n).$$

\no According to Theorem $4.2.1$ and the Euclid's fifth postulate,
we know that

$$\angle 2+\angle 4+\sum\limits_{i=2}^{n-1}f(x_i)= (n-2)\pi,$$

$$\angle 1+\angle 3\geq \pi$$

\no Therefore, we get that

$$\sum\limits_{i=1}^nf(x_{i})=(n-2)\pi+ \angle 1+\angle 3\geq (n-1)\pi.\quad\quad \natural$$

For opened $m$-lines with self-intersections, we know a result as
in the following.

\vskip 3mm

\no{\bf Theorem $4.2.3$} \ {Let $(M,\mu)$ be a planar map
geometry. An $m$-line $L$ in $(M,\mu)$ passing through $n$
non-euclidean points $x_1,x_2,\cdots, x_n$ is opened only with $l$
self-intersections if and only if there exist integers $i_j$ and
$s_{i_j}, 1\leq j\leq l$ with $1\leq i_j, s_{i,j}\leq n$ and
$i_j\not=i_t$ if $t\not=j$ such that}

$$(s_{i_j}-2)\pi \ < \ \sum\limits_{h=1}^{s_{i_j}}f(x_{i_j+h}) \ <  \ (s_{i_j}-1)\pi.$$

\vskip 3mm

{\it Proof} \ If an $m$-line $L$ passing through $m$-points
$x_{t+1},x_{t+2},\cdots, x_{t+s_t}$ only has one self-intersection
point, let us look at Fig.$4.4$ in where $x_{t+1}x_{t+s_t}$ is an
added auxiliary $m$-line segment.

\includegraphics[bb=10 10 400 88]{sgm58.eps}\vskip 2mm

\c{\bf Fig.$4.4$}

\vskip 2mm

\no We know that

$$\angle 1+\angle 2= f(x_{t+1}) \ {\rm and} \ \angle 3+\angle 4=f(x_{t+s_t}).$$

Similar to the proof of Theorem $4.2.2$, by Theorem $4.2.1$ and
the Euclid's fifth postulate, we know that

$$\angle 2+\angle 4+\sum\limits_{j=2}^{s_t-1}f(x_{t+j})= (s_t-2)\pi$$

\no and

$$\angle 1+\angle 3 \ < \ \pi.$$

\no Whence, we get that

$$(s_t-2)\pi \ < \ \sum\limits_{j=1}^{s_t}f(x_{t+j})\ < \ (s_t-1)\pi.$$

Therefore, if $L$ is opened only with $l$ self-intersection
points, we can find integers $i_j$ and $s_{i_j}, 1\leq j\leq l$
with $1\leq i_j, s_{i,j}\leq n$ and $i_j\not=i_t$ if $t\not=j$
such that $L$ passing through $x_{i_j+1}, x_{i_j+2},\cdots,
x_{i_j+s_j}$ only has one self-intersection point. By the previous
discussion, we know that

$$(s_{i_j}-2)\pi \ < \ \sum\limits_{h=1}^{s_{i_j}}f(x_{i_j+h}) \ <  \ (s_{i_j}-1)\pi.$$

\no This completes the proof. \quad\quad $\natural$

Notice that all $m$-lines considered in this section are consisted
by line segments or rays in an Euclid plane geometry. If the
length of each line segment tends to zero, then we get a curve at
the limitation in the usually sense. Whence, an $m$-line in a
planar map geometry can be also seen as a discretization for plane
curves and also has relation with differential equations. Readers
interested in those materials can see in Chapter $5$ for more
details.

\vskip 4mm

\no{\bf $4.2.2.$ Curvature of an $m$-line}

\vskip 3mm

\no The curvature at a point of a curve $C$ is a measure of how
quickly the tangent vector changes direction with respect to the
length of arc, such as those of the Gauss curvature, the Riemann
curvature, $\cdots$, etc.. In Fig.$4.5$ we present a smooth curve
and the changing of tangent vectors.

\includegraphics[bb=10 10 400 110]{sgm59.eps}

\vskip 3mm

\c{\bf Fig.$4.5$}

\vskip 2mm

To measure the changing of vector $v_1$ to $v_2$, a simpler way is
by the changing of the angle between vectors $v_1$ and $v_2$. If a
curve $C=f(s)$ is smooth, then the changing rate of the angle
between two tangent vector with respect to the length of arc,
i.e., $\frac{df}{ds}$ is continuous. For example, as we known in
the differential geometry, the Gauss curvature at every point of a
circle $x^2+y^2=r^2$ of radius $r$ is $\frac{1}{r}$. Whence, the
changing of the angle from vectors $v_1$ to $v_2$ is

$$\int\limits_{A}^B\frac{1}{r}ds=\frac{1}{r}|\widehat{AB}|=\frac{1}{r}r\theta =\theta.$$

\no By results in an Euclid plane geometry, we know that $\theta$
is also the angle between vectors $v_1$ and $v_2$. As we
illustrated in Subsection $4.2.1$, an $m$-line in a planar map
geometry is consisted by line segments or rays. Therefore, the
changing rate of the angle between two tangent vector with respect
to the length of arc is not continuous. Similar to the definition
of the set curvature in the reference $[1]$, we present a discrete
definition for the curvature of $m$-lines as follows.

\vskip 4mm

\no{\bf Definition $4.2.1$} \ {\it Let $L$ be an $m$-line in a
planar map geometry $(M,\mu)$ with the set $W$ of non-euclidean
points. The curvature $\omega(L)$ of $L$ is defined by}

$$\omega(L)=\sum\limits_{p\in W}(\pi - \varpi(p)),$$

\no{\it where $\varpi(p)=f(p)$ if $p$ is on an edge $(u,v)$ in map
$M$ on a plane $\sum$ with an angle function $f: \ \sum\rightarrow
\sum$.}

\vskip 3mm

In the classical differential geometry, the {\it Gauss mapping}
and the {\it Gauss curvature} on surfaces are defined as follows:

\vskip 3mm

{\it Let ${\mathcal S}\subset R^3$ be a surface with an
orientation {$\bf N$}. The mapping $ N: {\mathcal S}\rightarrow
S^2$ takes its value in the unit sphere}

$$S^2=\{(x,y,z)\in R^3|x^2+y^2+z^2=1\}$$

\no{\it along the orientation} {$\bf N$}. {\it The map $N:
{\mathcal S}\rightarrow S^2$, thus defined, is called a Gauss
mapping and the determinant of $K(p)=d{\bf N}_p$ a Gauss
curvature.}\vskip 2mm

We know that for a point $p\in {\mathcal S}$ such that the
Gaussian curvature $K(p)\not=0$ and a connected neighborhood $V$
of $p$ with $K$ does not change sign,

$$K(p)=\lim_{A\rightarrow 0}\frac{N(A)}{A},$$

\no where $A$ is the area of a region $B\subset V$ and $ N(A)$ is
the area of the image of $B$ by the Gauss mapping $ N: {\mathcal
S}\rightarrow S^2$.

The well-known {\it Gauss-Bonnet theorem} for a compact surface
says that

\vskip 3mm

$$\int\int_{\mathcal S}Kd\sigma =2\pi\chi (S),$$

\no{ for any orientable compact surface $S$.}

For a simply closed $m$-line, we also have a result similar to the
Gauss-Bonnet theorem, which can be also seen as a discrete
Gauss-Bonnet theorem on a plane.

\vskip 4mm

\no{\bf Theorem $4.2.4$} \ {\it Let $L$ be a simply closed
$m$-line passing through $n$ non-euclidean points $x_1,x_2,\cdots,
x_n$ in a planar map geometry $(M,\mu)$. Then}

$$\omega(L)=2\pi.$$

\vskip 3mm

{\it Proof} \ According to Theorem $4.2.1$, we know that

$$\sum\limits_{i=1}^nf(x_i)=(n-2)\pi,$$

\no where $f(x_i)$ denotes the angle function value at an
$m$-point $x_i, 1\leq i\leq n$. Whence, by Definition $4.2.1$ we
know that

\begin{eqnarray*}
\omega(L) &=& \sum\limits_{i=1}^n(\pi- f(x_i))\\
&=& \pi n-\sum\limits_{i=1}^nf(x_i)\\
&=& \pi n-(n-2)\pi = 2\pi.\quad\quad \natural
\end{eqnarray*}

Similarly, we get a result for the sum of curvatures on the planar
map $M$ in a planar geometry $(M,\mu)$.

\vskip 4mm

\no{\bf Theorem $4.2.6$} \ {\it Let $(M,\mu)$ be a planar map
geometry. Then the sum $\omega(M)$ of curvatures on edges in a map
$M$ is}

$$\omega(M)=2\pi s(M),$$

\no{\it where $s(M)$ denotes the sum of length of edges in $M$.}

\vskip 3mm

{\it Proof} \ Notice that the sum $\omega(u,v)$ of curvatures on
an edge $(u,v)$ of $M$ is

$$\omega(u,v)=\int\limits_{v}^u(\pi-f(s))ds=\pi|\widehat{(u,v)}|
-\int\limits_{v}^uf(s)ds.$$

Since $M$ is a planar map, each of its edges appears just two
times with an opposite direction. Whence, we get that

\begin{eqnarray*}
\omega(M)&=& \sum\limits_{(u,v)\in
E(M)}\omega(u,v)+\sum\limits_{(v,u)\in E(M)}\omega(v,u)\\
&=& \pi \sum\limits_{(u,v)\in
E(M)}(|\widehat{(u,v)}|+|\widehat{(v,u)}|)-(\int\limits_{v}^uf(s)ds+\int\limits_{u}^vf(s)ds)\\
&=& 2\pi s(M)\quad\quad \natural
\end{eqnarray*}

Notice that if we assume $s(M)=1$, then Theorem $4.2.6$ turns to
the Gauss-bonnet theorem for a sphere. Similarly, if we consider
general map geometry on an orientable surface, similar results can
be also obtained such as those materials in Problem $4.7.8$ and
Conjecture $4.7.1$ in the final section of this chapter.

\vskip 5mm

\no{\bf \S $4.3$ \ Polygons in a Planar Map Geometry}

\vskip 4mm

\no{\bf $4.3.1.$ Existence}\vskip 2mm

\no In an Euclid plane geometry, we have encountered triangles,
quadrilaterals, $\cdots$, and generally, $n$-polygons, i.e., these
graphs on a plane with $n$ straight line segments not on the same
line connected with one after another. There are no $1$ and
$2$-polygons in an Euclid plane geometry since every point is
euclidean. The definition of $n$-polygons in a planar map geometry
$(M,\mu)$ is similar to that of an Euclid plane geometry.

\vskip 4mm

\no{\bf Definition $4.3.1$} \ {\it An $n$-polygon in a planar map
geometry $(M,\mu)$ is defined to be a graph on  $(M,\mu)$ with $n$
$m$-line segments two by two without self-intersections and
connected with one after another.}

\vskip 3mm

Although their definition is similar, the situation is more
complex in a planar map geometry $(M,\mu)$. We have found a
necessary and sufficient condition for $1$-polygon in Theorem
$4.2.1$, i.e., $1$-polygons maybe exist in a planar map geometry.
In general, we can find $n$-polygons in a planar map geometry for
any integer $n, n\geq 1$.

Examples of polygon in a planar map geometry $(M,\mu)$ are shown
in Fig.$4.6$, in where $(a)$ is a $1$-polygon with $u,v,w$ and $t$
being non-euclidean points, $(b)$ is a $2$-polygon with vertices
$A,B$ and non-euclidean points $u,v$, $(c)$ is a triangle with
vertices $A,B,C$ and a non-euclidean point $u$ and $(d)$ is a
quadrilateral with vertices $A,B,C$ and $D$.

\includegraphics[bb=10 10 400 110]{sgm60.eps}

\vskip 3mm

\c{\bf Fig.$4.6$}

\vskip 4mm

\no{\bf Theorem $4.3.1$} \ {\it There exists a $1$-polygon in a
planar map geometry $(M,\mu)$ if and only if there are
non-euclidean points $u_1,u_2,\cdots,u_l$ with $l\geq 3$ such
that}

$$\sum\limits_{i=1}^lf(u_i)=(l-2)\pi,$$

\no{\it where $f(u_i)$ denotes the angle function value at the
point $u_i$, $1\leq i\leq l$.}

\vskip 3mm

{\it Proof} \ According to Theorem $4.2.1$, an $m$-line passing
through $l$ non-euclidean points $u_1,u_2,\cdots,u_l$ is simply
closed if and only if

$$\sum\limits_{i=1}^lf(u_i)=(l-2)\pi,$$

\no i.e., $1$-polygon exists in $(M,\mu)$ if and only if there are
non-euclidean points $u_1,u_2,\cdots,u_l$ with the above formula
hold.

Whence, we only need to prove $l\geq 3$. Since there are no
$1$-polygons or $2$-polygons in an Euclid plane geometry, we must
have $l\geq 3$ by the Hilbert's axiom $I-2$. In fact, for $l=3$ we
can really find a planar map geometry $(M,\mu)$ with a $1$-polygon
passing through three non-euclidean points $u, v$ and $w$. Look at
Fig.$4.7$,

\includegraphics[bb=10 10 400 110]{sgm61.eps}

\c{\bf Fig.$4.7$}

\vskip 2mm

\no in where the angle function values are
$f(u)=f(v)=f(w)=\frac{2}{3}\pi$ at $u, v$ and $w$.\quad\quad
$\natural$

Similarly, for $2$-polygons we get the following result.

\vskip 4mm

\no{\bf Theorem $4.3.2$} \ {\it There are $2$-polygons in a planar
map geometry $(M,\mu)$ only if there are at least one
non-euclidean point in  $(M,\mu)$.}

\vskip 3mm

{\it Proof} \ In fact, if there is a non-euclidean point $u$ in
$(M,\mu)$, then each straight line enter $u$ will turn an angle
$\theta =\pi-\frac{f(u)}{2}$ or $\frac{f(u)}{2}-\pi$ from the
initial straight line dependent on that $u$ is elliptic or
hyperbolic. Therefore, we can get a $2$-polygon in $(M,\mu)$ by
choice a straight line $AB$ passing through euclidean points in
$(M,\mu)$, such as the graph shown in Fig.$4.8$.

\includegraphics[bb=10 10 400 110]{sgm62.eps}

\c{\bf Fig.$4.8$}

\vskip 2mm

\no This completes the proof. \quad\quad $\natural$

For the existence of $n$-polygons with $n\geq 3$, we have a
general result as in the following.

\vskip 4mm

\no{\bf Theorem $4.3.3$} \ {\it For any integer $n, n\geq 3$,
there are $n$-polygons in a planar map geometry $(M,\mu)$.}

\vskip 3mm

{\it Proof} \ Since in an Euclid plane geometry, there are
$n$-polygons for any integer $n, n\geq 3$. Therefore, there are
also $n$-polygons in a planar map geometry $(M,\mu)$ for any
integer $n, n\geq 3$. \quad\quad $\natural$

\vskip 4mm

\no{\bf $4.3.2.$ Sum of internal angles}

\vskip 3mm

\no For the sum of the internal angles in an $n$-polygon, we have
the following result.

\vskip 4mm

\no{\bf Theorem $4.3.4$} \ {\it Let $\prod$ be an $n$-polygon in a
map geometry with its edges passing through non-euclidean points
$x_1,x_2,\cdots, x_l$. Then the sum of internal angles in $\prod$
is}

$$(n+l-2)\pi-\sum\limits_{i=1}^lf(x_i),$$

\no{\it where $f(x_i)$ denotes the value of the angle function $f$
at the point $x_i, 1\leq i\leq l$.}

\vskip 3mm

{\it Proof} \ Denote by $U, V$ the sets of elliptic points and
hyperbolic points in $x_1,x_2,\cdots, x_l$ and $|U|=p, |V|=q$,
respectively. If an $m$-line segment passes through an elliptic
point $u$, add an auxiliary line segment $AB$ in the plane as
shown in Fig.$4.9$(1). Then we get that

$$\angle{\rm a}= \angle 1 + \angle 2=\pi - f(u).$$

If an $m$-line passes through a hyperbolic point $v$, also add an
auxiliary line segment $AB$ in the plane as that shown in
Fig.$4.9$(2). Then we get that

$${\rm angle} \ b= {\rm angle} 3 + {\rm angle} 4=f(v)-\pi.$$

\vskip 3mm

\includegraphics[bb=5 5 100 130]{sgm63.eps}

\vskip 2mm

\c{\bf Fig.$4.9$}

Since the sum of internal angles of an $n$-polygon in a plane is
$(n-2)\pi$ whenever it is a convex or concave polygon, we know
that the sum of the internal angles in $\prod$ is

\begin{eqnarray*}
& \ &(n-2)\pi+\sum\limits_{x\in U}(\pi-f(x))-\sum\limits_{y\in V}(f(y)-\pi)\\
& \ & =(n+p+q-2)\pi-\sum\limits_{i=1}^lf(x_i)\\
& \ & =(n+l-2)\pi-\sum\limits_{i=1}^lf(x_i).
\end{eqnarray*}

\no This completes the proof. \quad\quad $\natural$

A triangle is called {\it euclidean, elliptic} or {\it hyperbolic}
if its edges only pass through one kind of euclidean, elliptic or
hyperbolic points. As a consequence of Theorem $4.3.4$, we get the
sum of the internal angles of a triangle in a map geometry which
is consistent with these already known results .

\vskip 4mm

\no{\bf Corollary $4.3.1$} {\it Let $\triangle$ be a triangle in a
planar map geometry $(M,\mu)$. Then }

($i$) {\it the sum of its internal angles is equal to $\pi$ if
$\triangle$ is euclidean;}

($ii$) {\it  the sum of its internal angles is less than $\pi$ if
$\triangle$ is elliptic;}

($iii$) {\it  the sum of its internal angles is more than $\pi$ if
$\triangle$ is hyperbolic.}

\vskip 3mm

{\it Proof} \ Notice that the sum of internal angles of a triangle
is

$$\pi+ \sum\limits_{i=1}^l(\pi-f(x_i))$$

\no if it passes through non-euclidean points $x_1,x_2,\cdots,
x_l$. By definition, if these $x_i, 1\leq i\leq l$ are one kind of
euclidean, elliptic, or hyperbolic, then we have that $f(x_i)=
\pi$, or $f(x_i)< \pi$, or $f(x_i)> \pi$ for any integer $i, 1\leq
i\leq l$. Whence, the sum of internal angles of an euclidean,
elliptic or hyperbolic triangle is equal to, or lees than, or more
than $\pi$.\quad\quad $\natural$

\vskip 4mm

\no{\bf $4.3.3.$ Area of a polygon}

\vskip 3mm

\no As it is well-known, calculation for the area $A(\triangle)$
of a triangle $\triangle$ with two sides $a,b$ and the value of
their include angle $\theta$ or three sides $a,b$ and $c$ in an
Euclid plane is simple. Formulae for its area are

$$ A(\triangle)=\frac{1}{2}ab\sin\theta \ {\rm or} \
A(\triangle)=\sqrt{s(s-a)(s-b)(s-c)},$$

\no where $s=\frac{1}{2}(a+b+c)$. But in a planar map geometry,
calculation for the area of a triangle is complex since each of
its edge maybe contains non-euclidean points. Where, we only
present a programming for calculation the area of a triangle in a
planar map geometry.

\vskip 3mm

{\bf STEP $1$} \ {\it Divide a triangle into triangles in an
Euclid plane such that no edges contain non-euclidean points
unless their endpoints;}\vskip 2mm

{\bf STEP $2$} \ {\it Calculate the area of each triangle;}\vskip
2mm

{\bf STEP $3$} \ {\it Sum up all of areas of these triangles to
get the area of the given triangle in a planar map geometry.}

\vskip 2mm

The simplest cases for triangle is the cases with only one
non-euclidean point such as those shown in Fig.$4.10(1)$ and $(2)$
with an elliptic point $u$ or with a hyperbolic point $v$.

\includegraphics[bb=5 5 100 105]{sgm64.eps}

\c{\bf Fig.$4.10$}

\vskip 2mm

\no Add an auxiliary line segment $AB$ in Fig.$4.10$. Then by
formulae in the plane trigonometry, we know that

$$A(\triangle
ABC)=\sqrt{s_1(s_1-a)(s_1-b)(s_1-t)}+\sqrt{s_2(s_2-c)(s_2-d)(s_2-t)}$$

\no for case $(1)$ in Fig.$4.10$ and

$$A(\triangle
ABC)=\sqrt{s_1(s_1-a)(s_1-b)(s_1-t)}-\sqrt{s_2(s_2-c)(s_2-d)(s_2-t)}$$

\no for case $(2)$ in Fig.$4.10$, where

$$t=\sqrt{c^2+d^2-2cd\cos\frac{f(x)}{2}}$$

\no with $x=u$ or $v$ and

$$s_1=\frac{1}{2}(a+b+t), \ \ s_2=\frac{1}{2}(c+d+t).$$

Generally, let $\triangle ABC$ be a triangle with its edge $AB$
passing through $p$ elliptic or $p$ hyperbolic points
$x_1,x_2,\cdots, x_p$ simultaneously, as those shown in
Fig.$4.11(1)$ and $(2)$.

\includegraphics[bb=5 5 100 140]{sgm65.eps}

\c{\bf Fig.$4.11$}\vskip 2mm

\no Where $|AC|=a$, $|BC|=b$ and $|Ax_1|=c_1$,
$|x_1x_2|=c_2,\cdots, |x_{p-1}x_p|=c_p$ and $|x_pB|=c_{p+1}$.
Adding auxiliary line segments $Ax_2, Ax_3,\cdots, Ax_p,AB$ in
Fig.$4.11$, then we can find its area by the programming STEP $1$
to STEP $3$. By formulae in the plane trigonometry, we get that

$$|Ax_2|=\sqrt{c_1^2+c_2^2-2c_1c_2\cos\frac{f(x_1)}{2}},$$

$$\angle Ax_2x_1= \cos^{-1}\frac{c_1^2-c_2^2-|Ax_1|^2}{2c_2|Ax_2|},$$

$$\angle Ax_2x_3= \frac{f(x_2)}{2}-\angle Ax_2x_1 \ \ {\rm or}  \ \
2\pi-\frac{f(x_2)}{2}-\angle Ax_2x_1,$$

$$|Ax_3|=\sqrt{|Ax_2|^2+c_3^2-2|Ax_2|c_3\cos(\frac{f(x_2)}{2}-\angle Ax_2x_3)},$$

$$\angle Ax_3x_2= \cos^{-1}\frac{|Ax_2|^2-c_3^2-|Ax_3|^2}{2c_3|Ax_3|},$$

$$\angle Ax_2x_3= \frac{f(x_3)}{2}-\angle Ax_3x_2\ \ {\rm or}  \ \
2\pi-\frac{f(x_3)}{2}-\angle Ax_3x_2,$$

$$\cdots\cdots\cdots\cdots\cdots\cdots\cdots\cdots\cdots$$

\no and generally, we get that

$$|AB|=\sqrt{|Ax_p|^2+c_{p+1}^2-2|Ax_p|c_{p+1}\cos\angle Ax_pB}.$$

\no Then the area of the triangle $\triangle ABC$ is

\begin{eqnarray*}
A(\triangle ABC)&=& \sqrt{s_p(s_p-a)(s_p-b)(s_p-|AB|)}\\
&+&\sum\limits_{i=1}^p\sqrt{s_i(s_i-|Ax_i|)(s_i-c_{i+1})(s_i-|Ax_{i+1}|)}
\end{eqnarray*}

\no for case $(1)$ in Fig.$4.11$ and

\begin{eqnarray*}
A(\triangle ABC)&=& \sqrt{s_p(s_p-a)(s_p-b)(s_p-|AB|)}\\
&-&
\sum\limits_{i=1}^p\sqrt{s_i(s_i-|Ax_i|)(s_i-c_{i+1})(s_i-|Ax_{i+1}|)}
\end{eqnarray*}

\no for case $(2)$ in Fig.$4.11$, where for any integer $i,1\leq
i\leq p-1$,

$$s_i=\frac{1}{2}(|Ax_i|+c_{i+1}+|Ax_{i+1}|)$$

\no and

$$s_p=\frac{1}{2}(a+b+|AB|).$$

Certainly, this programming can be also applied to calculate the
area of an $n$-polygon in a planar map geometry in general.

\vskip 5mm

\no{\bf \S $4.4$ \ Circles in a Planar Map Geometry}

\vskip 4mm

\no The length of an $m$-line segment in a planar map geometry is
defined in the following definition.

\vskip 4mm

\no{\bf Definition $4.4.1$} \ {\it The length $|AB|$ of an
$m$-line segment $AB$ consisted by $k$ straight line segments
$AC_1,C_1C_2$, $C_2C_3,\cdots$,$C_{k-1}B$ in a planar map geometry
$(M,\mu)$ is defined by}

$$|AB|= |AC_1|+|C_1C_2|+|C_2C_3|+\cdots+|C_{k-1}B|.$$

As that shown in Chapter $3$, there are not always exist a circle
with any center and a given radius in a planar map geometry in the
sense of the Euclid's definition. Since we have introduced angle
function on a planar map geometry, we can likewise the Euclid's
definition to define an $m$-circle in a planar map geometry in the
next definition.

\vskip 4mm

\no{\bf Definition $4.4.2$} \ {\it A closed curve $C$ without
self-intersection in a planar map geometry $(M,\mu)$ is called an
$m$-circle if there exists an $m$-point $O$ in $(M,\mu)$ and a
real number $r$ such that $|OP|=r$ for each $m$-point $P$ on $C$.}

\vskip 3mm

Two Examples for $m$-circles in a planar map geometry $(M,\mu)$
are shown in Fig.$4.12(1)$ and $(2)$. The $m$-circle in
Fig.$4.12(1)$ is a circle in the Euclid's sense, but $(2)$ is not.
Notice that in Fig.$4.12(2)$, $m$-points $u$ and $v$ are elliptic
and the length $|OQ|=|Ou|+|uQ|=r$ for an $m$-point $Q$ on the
$m$-circle $C$, which seems likely an ellipse but it is not. The
$m$-circle $C$ in Fig.$4.12(2)$ also implied that $m$-circles are
more complex than those in an Euclid plane geometry.

\includegraphics[bb=5 5 100 120]{sgm66.eps}

\vskip 2mm

\c{\bf Fig.$4.12$}\vskip 2mm

We have a necessary and sufficient condition for the existence of
an $m$-circle in a planar map geometry.

\vskip 4mm

\no{\bf Theorem $4.4.1$} \ {\it Let $(M,\mu)$ be a planar map
geometry on a plane $\sum$ and $O$ an $m$-point on $(M,\mu)$. For
a real number $r$, there is an $m$-circle of radius $r$ with
center $O$ if and only if $O$ is in a non-outer face of $M$ or $O$
is in the outer face of $M$ but for any $\epsilon, r> \epsilon
>0$, the initial and final intersection points of a circle of
radius $\epsilon$ with $M$ in an Euclid plane $\sum$ are euclidean
points.}

\vskip 3mm

{\it Proof} \ If there is a solitary non-euclidean point $A$ with
$|OA|< r$, then by those materials in Chapter $3$, there are no
$m$-circles in $(M,\mu)$ of radius $r$ with center $O$.

Now if $O$ is in the outer face of $M$ but there exists a number
$\epsilon, r >\epsilon > 0$ such that one of the initial and final
intersection points of a circle of radius $\epsilon$ with $M$ on
$\sum$ is non-euclidean point, then points with distance $r$ to
$O$ in $(M,\mu)$ at least has a gap in a circle with an Euclid
sense. See Fig.$4.13$ for details, in where $u$ is a non-euclidean
point and the shade field denotes the map $M$. Therefore there are
no $m$-circles in $(M,\mu)$ of radius $r$ with center $O$.

\includegraphics[bb=5 5 100 120]{sgm67.eps}

\vskip 2mm

\c{\bf Fig.$4.13$}\vskip 2mm

Now if $O$ in the outer face of $M$ but for any $\epsilon, r>
\epsilon >0$, the initial and final intersection points of a
circle of radius $\epsilon$ with $M$ on $\sum$ are euclidean
points or $O$ is in a non-outer face of $M$, then by the
definition of angle functions, we know that all points with
distance $r$ to $O$ is a closed smooth curve on $\sum$, for
example, see Fig.$4.14(1)$ and $(2)$.

\includegraphics[bb=5 5 100 120]{sgm68.eps}

\vskip 2mm

\c{\bf Fig.$4.14$}\vskip 2mm

\no Whence it is an $m$-circle. \quad\quad $\natural$

We construct a polar axis $OX$ with center $O$ in a planar map
geometry as that in an Euclid geometry. Then each $m$-point $A$
has a coordinate $(\rho,\theta)$, where $\rho$ is the length of
the $m$-line segment $OA$ and $\theta$ is the angle between $OX$
and the straight line segment of $OA$ containing the point $A$. We
get an equation for an $m$-circle of radius $r$ which has the same
form as that in the analytic geometry of plane.

\vskip 4mm

\no{\bf Theorem $4.4.2$} \ {\it In a planar geometry $(M,\mu)$
with a polar axis $OX$ of center $O$, the equation of each
$m$-circle of radius $r$ with center $O$ is}

$$\rho = r.$$

\vskip 3mm

{\it Proof} \ By the definition of an $m$-circle $C$ of radius
$r$, every $m$-point on $C$ has a distance $r$ to its center $O$.
Whence, its equation is $\rho = r$ in a planar map geometry with a
polar axis $OX$ of center $O$.\quad\quad $\natural$

\vskip 5mm

\no{\bf \S $4.5$ \ Line Bundles in a Planar Map Geometry}

\vskip 4mm

\no The behaviors of $m$-line bundles is need to clarify from a
geometrical sense. Among those $m$-line bundles the most important
is parallel bundles defined in the next definition, which is also
motivated by the Euclid's fifth postulate discussed in the
reference $[54]$ first.

\vskip 4mm

\no{\bf Definition $4.5.1$} \ {\it A family $\mathcal L$ of
infinite $m$-lines not intersecting each other in a planar
geometry $(M,\mu)$ is called a parallel bundle.}

\vskip 2mm

In Fig.$4.15$, we present all cases of parallel bundles passing
through an edge in planar geometries, where, (a) is the case with
the same type points $u,v$ and $\rho_M (u)\mu (u)=\rho_M (v)\mu
(v)=2\pi$, (b) and (c) are the same type cases with $\rho_M (u)\mu
(u)
> \rho_M (v)\mu (v)$ or $\rho_M (u)\mu (u)=\rho_M (v)\mu
(v)>2\pi$ or $< 2\pi$ and (d) is the case with an elliptic point
$u$ but a hyperbolic point $v$.

\vskip 3mm

\includegraphics[bb=5 5 100 130]{sgm69.eps}

\vskip 2mm

\c{\bf Fig.$4.15$}

\no Here, we assume the angle at the intersection point is in
clockwise, that is, a line passing through an elliptic point will
bend up and passing through a hyperbolic point will bend down,
such as those cases (b),(c) in the Fig.$4.15$. Generally, we
define a {\it sign function $sign(f)$ of an angle function $f$} as
follows.

\vskip 4mm

\no{\bf Definition $4.5.2$} \ {\it For a vector
$\overrightarrow{O}$ on the Euclid plane called an orientation, a
sign function $sign(f)$ of an angle function $f$ at an $m$-point
$u$ is defined by}

\[
sign(f)(u)=\left\{
\begin{array}{lr}
1, & {\rm if} \ u \ {\rm is \ elliptic,}\\
0, & {\rm if} \ u \ {\rm is \ euclidean,}\\
-1, & {\rm if} \ u \ {\rm is \ hyperbolic}.
\end{array}
\right.
\]

\vskip 3mm

We classify parallel bundles in planar map geometries along an
orientation $\overrightarrow{O}$ in this section.

\vskip 4mm

\no{\bf $4.5.1.$ A condition for parallel bundles}

\vskip 3mm

\no We investigate the behaviors of parallel bundles in a planar
map geometry $(M,\mu)$. Denote by $f(x)$ the angle function value
at an intersection $m$-point of an $m$-line $L$ with an edge
$(u,v)$ of $M$ and a distance $x$ to $u$ on $(u,v)$ as shown in
Fig.$4.15(a)$. Then we get an elementary result as in the
following.

\vskip 4mm

\no{\bf Theorem $4.5.1$} \ {\it A family $\mathcal L$ of parallel
$m$-lines passing through an edge $(u,v)$ is a parallel bundle if
and only if }

$$\left. \frac{df}{dx} \right|_+\geq 0.$$

\vskip 3mm

{\it Proof} \ If $\mathcal L$ is a parallel bundle, then any two
$m$-lines $L_1,L_2$ will not intersect after them passing through
the edge $uv$. Therefore, if $\theta_1, \theta_2$ are the angles
of $L_1, L_2$ at the intersection $m$-points of $L_1, L_2$ with
$(u,v)$ and $L_2$ is far from $u$ than $L_1$, then we know
$\theta_2\geq \theta_1$. Thereby we know that

$$f(x+\Delta x)-f(x)\geq 0$$

\no for any point with distance $x$ from $u$ and $\Delta x \ > \
0$. Therefore, we get that

$$\left. \frac{df}{dx} \right|_+=\lim_{\Delta x\to +0}\frac{f(x+\Delta x)-f(x)}
{\Delta x}\geq 0.$$

\no As that shown in the Fig.$4.15$.

Now if $\left.\frac{df}{dx}\right|_+\geq 0$, then $f(y)\geq f(x)$
if $y\geq x$. Since ${\mathcal L}$ is a family of parallel
$m$-lines before meeting $uv$, any two $m$-lines in ${\mathcal L}$
will not intersect each other after them passing through $(u,v)$.
Therefore, ${\mathcal L}$ is a parallel bundle.\quad\quad
$\natural$

A general condition for a family of parallel $m$-lines passing
through a cut of a planar map being a parallel bundle is the
following.

\vskip 4mm

\no{\bf Theorem $4.5.2$} \ {\it Let $(M,\mu)$ be a planar map
geometry, $C=\{(u_1,v_1),(u_2,v_2),\cdots,$ $(u_l,v_l)\}$ a cut of
the map $M$ with order $(u_1,v_1),(u_2,v_2),\cdots,(u_l,v_l)$ from
the left to the right, $l\geq 1$ and the angle functions on them
are $f_1, f_2,\cdots ,f_l$ (also seeing Fig.$4.16$), respectively.

\includegraphics[bb=5 5 100 130]{sgm70.eps}

\vskip 2mm

\c{\bf Fig.$4.16$}

\no Then a family $\mathcal L$ of parallel $m$-lines passing
through $C$ is a parallel bundle if and only if for any $x, x \geq
0$,}

\begin{eqnarray*}
sign(f_1)(x)f'_{1+}(x)\geq 0\\
sign(f_1)(x)f'_{1+}(x)+sign(f_2)(x)f'_{2+}(x)\geq 0\\
sign(f_1)(x)f'_{1+}(x)+sign(f_2)(x)f'_{2+}(x)+sign(f_3)(x)f'_{3+}(x)\geq 0\\
\cdots\cdots\cdots\cdots\\
sign(f_1)(x)f'_{1+}(x)+sign(f_2)(x)f'_{2+}(x)+\cdots
+sign(f_l)(x)f'_{l+}(x)\geq 0.
\end{eqnarray*}

\vskip 2mm

{\it Proof} \ According to Theorem $4.5.1$, we know that $m$-lines
will not intersect after them passing through $(u_1,v_1)$ and
$(u_2,v_2)$  if and only if for $\forall \Delta x > 0$ and $x \geq
0$,

$$sign(f_2)(x)f_2(x+\Delta x)+sign(f_1)(x)f'_{1+}(x)\Delta x\geq sign(f_2)(x)f_2(x),$$

\no seeing Fig.$4.17$ for an explanation.

\includegraphics[bb=5 5 100 130]{sgm71.eps}

\vskip 2mm

\c{\bf Fig.$4.17$}\vskip 2mm

\no That is,

$$sign(f_1)(x)f'_{1+}(x)+sign(f_2)(x)f'_{2+}(x)\geq 0.$$

Similarly, $m$-lines will not intersect after them passing through
$(u_1,v_1), (u_2,v_2)$ and $(u_3,v_3)$ if and only if for $\forall
\Delta x
> 0$ and $x \geq 0$,

\begin{eqnarray*}
& \ & sign(f_3)(x)f_3(x+\Delta x)+sign(f_2)(x)f'_{2+}(x)\Delta
x\\
& \ & +sign(f_1)(x)f'_{1+}(x) \Delta x\geq sign(f_3)(x)f_3(x).
\end{eqnarray*}

\no That is,

$$sign(f_1)(x)f'_{1+}(x)+sign(f_2)(x)f'_{2+}(x)+sign(f_3)(x)f'_{3+}(x)\geq 0.$$

\no Generally, $m$-lines will not intersect after them passing
through $(u_1,v_1),(u_2,v_2),\cdots ,$ $(u_{l-1},v_{l-1})$ and
$(u_l,v_l)$ if and only if for $\forall \Delta x > 0$ and $x \geq
0$,

\begin{eqnarray*}
& \ & sign(f_l)(x)f_l(x+\Delta
x)+sign(f_{l-1})(x)f'_{l-1+}(x)\Delta x+\\
& \ & \cdots +sign(f_1)(x)f'_{1+}(x)\Delta x \geq
sign(f_l)(x)f_l(x).\end{eqnarray*}

\no Whence, we get that

$$sign(f_1)(x)f'_{1+}(x)+sign(f_2)(x)f'_{2+}(x)+\cdots +sign(f_l)(x)f'_{l+}(x)\geq 0.$$

Therefore, a family $\mathcal L$ of parallel $m$-lines passing
through $C$ is a parallel bundle if and only if for any $x, x \geq
0$, we have that

\begin{eqnarray*}
sign(f_1)(x)f'_{1+}(x)\geq 0\\
sign(f_1)(x)f'_{1+}(x)+sign(f_2)(x)f'_{2+}(x)\geq 0\\
sign(f_1)(x)f'_{1+}(x)+sign(f_2)(x)f'_{2+}(x)+sign(f_1)(x)f'_{3+}(x)\geq 0\\
\cdots\cdots\cdots\cdots\\
sign(f_1)(x)f'_{1+}(x)+sign(f_2)(x)f'_{2+}(x)+\cdots
+sign(f_1)(x)f'_{l+}(x)\geq 0.
\end{eqnarray*}

\no This completes the proof. \quad\quad $\natural$.

\vskip 4mm

\no{\bf Corollary $4.5.1$} \ {\it Let $(M,\mu)$ be a planar map
geometry, $C=\{(u_1,v_1),(u_2,v_2),\cdots ,$ $(u_l,v_l)\}$ a cut
of the map $M$ with order $(u_1,v_1),(u_2,v_2),\cdots ,(u_l,v_l)$
from the left to the right, $l\geq 1$ and the angle functions on
them are $f_1, f_2,\cdots ,f_l$, respectively. Then a family
$\mathcal L$ of parallel lines passing through $C$ is still
parallel lines after them leaving $C$ if and only if for any $x, x
\geq 0$,}

\begin{eqnarray*}
sign(f_1)(x)f'_{1+}(x)\geq 0\\
sign(f_1)(x)f'_{1+}(x)+sign(f_2)(x)f'_{2+}(x)\geq 0\\
sign(f_1)(x)f'_{1+}(x)+sign(f_2)(x)f'_{2+}(x)+sign(f_1)(x)f'_{3+}(x)\geq 0\\
\cdots\cdots\cdots\cdots\\
sign(f_1)(x)f'_{1+}(x)+sign(f_2)(x)f'_{2+}(x)+\cdots
+sign(f_1)(x)f'_{l-1+}(x)\geq 0.
\end{eqnarray*}

\no{\it and}

$$sign(f_1)(x)f'_{1+}(x)+sign(f_2)(x)f'_{2+}(x)+\cdots
+sign(f_1)(x)f'_{l+}(x)=0.$$

\vskip 3mm

{\it Proof} \ According to Theorem $4.5.2$, we know the condition
is a necessary and sufficient condition for $\mathcal L$ being a
parallel bundle. Now since lines in $\mathcal L$ are parallel
lines after them leaving $C$ if and only if for any $x\geq 0$ and
$\Delta x\geq 0$, there must be that

$$sign(f_l)f_l(x+\Delta x)+sign(f_{l-1})f'_{l-1+}(x)\Delta x+\cdots +sign(f_1)f'_{1+}(x)\Delta x
= sign(f_l)f_l(x).$$

\no Therefore, we get that

$$sign(f_1)(x)f'_{1+}(x)+sign(f_2)(x)f'_{2+}(x)+\cdots
+sign(f_1)(x)f'_{l+}(x)=0.\quad\quad \natural$$

When do some parallel $m$-lines parallel the initial parallel
lines after them passing through a cut $C$ in a planar map
geometry? The answer is in the next result.

\vskip 4mm

\no{\bf Theorem $4.5.3$} \ {\it Let $(M,\mu)$ be a planar map
geometry, $C=\{(u_1,v_1),(u_2,v_2),\cdots ,$ $(u_l,v_l)\}$ a cut
of the map $M$ with order $(u_1,v_1),(u_2,v_2),\cdots ,(u_l,v_l)$
from the left to the right, $l\geq 1$ and the angle functions on
them are $f_1, f_2,\cdots ,f_l$, respectively. Then the parallel
$m$-lines parallel the initial parallel lines after them passing
through $C$ if and only if for $\forall x\geq 0$,}

\begin{eqnarray*}
sign(f_1)(x)f'_{1+}(x)\geq 0\\
sign(f_1)(x)f'_{1+}(x)+sign(f_2)(x)f'_{2+}(x)\geq 0\\
sign(f_1)(x)f'_{1+}(x)+sign(f_2)(x)f'_{2+}(x)+sign(f_1)(x)f'_{3+}(x)\geq 0\\
\cdots\cdots\cdots\cdots\\
sign(f_1)(x)f'_{1+}(x)+sign(f_2)(x)f'_{2+}(x)+\cdots
+sign(f_1)(x)f'_{l-1+}(x)\geq 0.
\end{eqnarray*}

\no{\it and}

$$sign(f_1)f_1(x)+sign(f_2)f_2(x)+\cdots +sign(f_1)(x)f_l(x)=l\pi.$$

\vskip 2mm

{\it Proof} \ According to Theorem $4.5.2$ and Corollary $4.5.1$,
we know that these parallel $m$-lines satisfying conditions of
this theorem is a parallel bundle.

We calculate the angle $\alpha (i,x)$ of an $m$-line $L$ passing
through an edge $u_iv_i, 1\leq i\leq l$ with the line before it
meeting $C$ at the intersection of $L$ with the edge $(u_i,v_i)$,
where $x$ is the distance of the intersection point to $u_1$ on
$(u_1,v_1)$, see also Fig.$4.18$. By definition, we know the angle
$\alpha (1,x)=sign(f_1)f(x)$ and $\alpha
(2,x)=sign(f_2)f_2(x)-(\pi-sign(f_1)f_1(x))=sign(f_1)f_1(x)+sign(f_2)f_2(x)-\pi$.

Now if $\alpha (i,x)=sign(f_1)f_1(x)+sign(f_2)f_2(x)+\cdots
+sign(f_i)f_i(x)-(i-1)\pi$, then we know that $\alpha
(i+1,x)=sign(f_{i+1})f_{i+1}(x)-(\pi -\alpha
(i,x))=sign(f_{i+1})f_{i+1}(x)+\alpha (i,x)-\pi$ similar to the
case $i=2$. Thereby we get that

$$\alpha (i+1,x)=sign(f_1)f_1(x)+sign(f_2)f_2(x)+\cdots +sign(f_{i+1})f_{i+1}(x)-i\pi.$$

Notice that an $m$-line $L$ parallel the initial parallel line
after it passing through $C$ if and only if $\alpha (l,x)=\pi$,
i.e.,

$$sign(f_1)f_1(x)+sign(f_2)f_2(x)+\cdots +sign(f_l)f_l(x)=l\pi.$$

\no This completes the proof. \quad\quad $\natural$

\vskip 4mm

\no{\bf $4.5.2.$ Linear conditions and combinatorial realization
for parallel bundles}

\vskip 3mm

\no For the simplicity, we can assume even that the function
$f(x)$ is linear and denoted it by $f_l(x)$. We calculate $f_l(x)$
in the first.

\vskip 4mm

\no{\bf Theorem $4.5.4$} \ {\it The angle function $f_l(x)$ of  an
$m$-line $L$ passing through an edge $(u,v)$ at a point with
distance $x$ to u is }

 $$f_l(x)=
(1-\frac{x}{d(u,v)})\frac{\rho (u)\mu
(v)}{2}+\frac{x}{d(u,v)}\frac{\rho (v)\mu (v)}{2},$$

\no{\it where, $d(u,v)$ is the length of the edge $(u,v)$.}

\vskip 3mm

{\it Proof} \ Since $f_l(x)$ is linear, we know that $f_l(x)$
satisfies the following equation.

$$\frac{f_l(x)-\frac{\rho (u)\mu (u)}{2}}{\frac{\rho (v)\mu (v)}{2}-
\frac{\rho (u)\mu (u)}{2}}=\frac{x}{d(u,v)},$$

\no Calculation shows that

 $$f_l(x)=
(1-\frac{x}{d(u,v)})\frac{\rho (u)\mu
(v)}{2}+\frac{x}{d(u,v)}\frac{\rho (v)\mu (v)}{2}.\quad \natural$$

\vskip 4mm

\no{\bf Corollary $4.5.2$} \ {\it Under the linear assumption, a
family $\mathcal L$ of parallel $m$-lines passing through an edge
$(u,v)$ is a parallel bundle if and only if}

$$\frac{\rho (u)}{\rho (v)}\leq\frac{\mu (v)}{\mu (u)}.$$

\vskip 3mm

{\it Proof} \ According to Theorem $4.5.1$, a family of parallel
$m$-lines passing through an edge $(u,v)$ is a parallel bundle if
and only if $f'(x)\geq 0$ for $\forall x, x \geq 0$, i.e.,

$$\frac{\rho (v)\mu (v)}{2d(u,v)}-\frac{\rho (u)\mu (u)}{2d(u,v)}\geq 0.$$

Therefore, a family $\mathcal L$ of parallel $m$-lines passing
through an edge $(u,v)$ is a parallel bundle if and only if

$$ \rho (v)\mu (v)\geq \rho (u)\mu (u).$$

\no Whence,

$$\frac{\rho (u)}{\rho (v)}\leq\frac{\mu (v)}{\mu (u)}.\quad\quad \natural$$

For a family of parallel $m$-lines passing through a cut, we get
the following condition.

\vskip 4mm

\no{\bf Theorem $4.5.5$} \ {\it Let $(M,\mu)$ be a planar map
geometry, $C=\{(u_1,v_1),(u_2,v_2),\cdots ,$ $(u_l,v_l)\}$ a cut
of the map $M$ with order $(u_1,v_1),(u_2,v_2),\cdots ,(u_l,v_l)$
from the left to the right, $l\geq 1$. Then under the linear
assumption, a family $L$ of parallel $m$-lines passing through $C$
is a parallel bundle if and only if the angle factor $\mu$
satisfies the following linear inequality system}

$$\rho (v_1)\mu (v_1)\geq \rho (u_1)\mu (u_1)$$

$$\frac{\rho (v_1)\mu (v_1)}{d(u_1,v_1)}+\frac{\rho (v_2)\mu
(v_2)}{d(u_2,v_2)}\geq \frac{\rho (u_1)\mu
(u_1)}{d(u_1,v_1)}+\frac{\rho (u_2)\mu (u_2)}{d(u_2,v_2)}$$

$$\cdots\cdots\cdots\cdots$$

\begin{eqnarray*}
\frac{\rho (v_1)\mu (v_1)}{d(u_1,v_1)}&+&\frac{\rho (v_2)\mu
(v_2)}{d(u_2,v_2)}+\cdots +\frac{\rho(v_l)\mu (v_l)}{d(u_l,v_l)}\\
&\geq &\frac{\rho (u_1)\mu (u_1)}{d(u_1,v_1)}+\frac{\rho (u_2)\mu
(u_2)}{d(u_2,v_2)}+\cdots + \frac{\rho (u_l)\mu
(u_l)}{d(u_l,v_l)}.
\end{eqnarray*}

\vskip 3mm

{\it Proof} \ Under the linear assumption, for any integer $i,
i\geq 1$ we know that

$$f'_{i+}(x)=\frac{\rho (v_i)\mu (v_i)-\rho (u_i)\mu (u_i)}{2d(u_i,v_i)}$$

\no by Theorem $4.5.4$. Thereby according to Theorem $4.5.2$, we
get that a family $L$ of parallel $m$-lines passing through $C$ is
a parallel bundle if and only if the angle factor $\mu$ satisfies
the following linear inequality system

$$\rho (v_1)\mu (v_1)\geq \rho (u_1)\mu (u_1)$$

$$\frac{\rho (v_1)\mu (v_1)}{d(u_1,v_1)}+\frac{\rho (v_2)\mu
(v_2)}{d(u_2,v_2)}\geq \frac{\rho (u_1)\mu
(u_1)}{d(u_1,v_1)}+\frac{\rho (u_2)\mu (u_2)}{d(u_2,v_2)}$$

$$\cdots\cdots\cdots\cdots$$

\begin{eqnarray*}
\frac{\rho (v_1)\mu (v_1)}{d(u_1,v_1)}&+&\frac{\rho (v_2)\mu
(v_2)}{d(u_2,v_2)}+\cdots +\frac{\rho(v_l)\mu (v_l)}{d(u_l,v_l)}\\
&\geq &\frac{\rho (u_1)\mu (u_1)}{d(u_1,v_1)}+\frac{\rho (u_2)\mu
(u_2)}{d(u_2,v_2)}+\cdots + \frac{\rho (u_l)\mu
(u_l)}{d(u_l,v_l)}.\end{eqnarray*}

\no This completes the proof. \quad\quad $\natural$

For planar maps underlying a regular graph, we have an interesting
consequence for parallel bundles in the following.

\vskip4mm

\no{\bf Corollary $4.5.3$} \ {\it Let $(M,\mu)$ be a planar map
geometry with $M$ underlying a regular graph,
$C=\{(u_1,v_1),(u_2,v_2),\cdots ,(u_l,v_l)\}$ a cut of the map $M$
with order $(u_1,v_1),(u_2,v_2),\cdots ,(u_l,v_l)$ from the left
to the right, $l\geq 1$. Then under the linear assumption, a
family $L$ of parallel lines passing through $C$ is a parallel
bundle if and only if the angle factor $\mu$ satisfies the
following linear inequality system.}

$$\mu (v_1)\geq  \mu (u_1)$$

$$\frac{\mu (v_1)}{d(u_1,v_1)}+\frac{\mu (v_2)}{d(u_2,v_2)}\geq
\frac{\mu (u_1)}{d(u_1,v_1)}+\frac{\mu (u_2)}{d(u_2,v_2)}$$

$$\cdots\cdots\cdots\cdots$$

$$ \frac{\mu (v_1)}{d(u_1,v_1)}+\frac{\mu (v_2)}{d(u_2,v_2)}+\cdots
+\frac{\mu (v_l)}{d(u_l,v_l)}\geq \frac{\mu
(u_1)}{d(u_1,v_1)}+\frac{\mu (u_2)}{d(u_2,v_2)}+\cdots + \frac{\mu
(u_l)}{d(u_l,v_l)}$$

\no{\it and particularly, if assume that all the lengths of edges
in $C$ are the same, then }

\begin{eqnarray*}
\mu (v_1)&\geq& \mu (u_1)\\
\mu (v_1)+\mu (v_2)&\geq& \mu (u_1)+\mu (u_2)\\
\cdots\cdots &\cdots& \cdots\cdots\\
\mu (v_1)+\mu (v_2)+\cdots +\mu (v_l)&\geq& \mu (u_1)+\mu
(u_2)+\cdots +\mu (u_l).
\end{eqnarray*}

\vskip 3mm

Certainly, by choice different angle factors, we can also get
combinatorial conditions for the existence of parallel bundles
under the linear assumption.

\vskip 4mm

\no{\bf Theorem $4.5.6$} \ {\it Let $(M,\mu)$ be a planar map
geometry, $C=\{(u_1,v_1),(u_2,v_2),\cdots ,$ $(u_l,v_l)\}$ a cut
of the map $M$ with order $(u_1,v_1),(u_2,v_2),\cdots ,(u_l,v_l)$
from the left to the right, $l\geq 1$. If}

$$\frac{\rho (u_i)}{\rho (v_i)}\leq\frac{\mu (v_i)}{\mu (u_i)}$$

\no{\it for any integer $i,i\geq 1$, then a family $L$ of parallel
$m$-lines passing through $C$ is a parallel bundle under the
linear assumption.}

\vskip 3mm

{\it Proof} \ Under the linear assumption we know that

$$f'_{i+}(x)=\frac{\rho (v_i)\mu (v_i)-\rho (u_i)\mu (u_i)}{2d(u_i,v_i)}$$

\no for any integer $i, i\geq 1$ by Theorem $4.5.4$. Thereby
$f'_{i+}(x)\geq 0$ for $i=1,2,\cdots ,l$. We get that

\begin{eqnarray*}
f'_1(x)\geq 0\\
f'_{1+}(x)+f'_{2+}(x)\geq 0\\
f'_{1+}(x)+f'_{2+}(x)+f'_{3+}(x)\geq 0\\
\cdots\cdots\cdots\cdots\\
f'_{1+}(x)+f'_{2+}(x)+\cdots +f'_{l+}(x)\geq 0.
\end{eqnarray*}

By Theorem $4.5.2$ we know that a family $L$ of parallel $m$-lines
passing through $C$ is still a parallel bundle.\quad\quad
$\natural$

\vskip 5mm

\no{\bf \S $4.6$ \ Examples of Planar Map Geometries}

\vskip 4mm

\no By choice different planar maps and define angle factors on
their vertices, we can get various planar map geometries. In this
section, we present some concrete examples for planar map
geometries.

\vskip 4mm

\no{\bf Example $4.6.1$} \ {\it A complete planar map $K_4$}

\vskip 3mm

We take a complete map $K_4$ embedded on a plane $\sum$ with
vertices $u,v,w$ and $t$ and angle factors

$$\mu(u)=\frac{\pi}{2}, \ \ \mu(v)=\mu(w)=\pi  \ {\rm and} \ \mu(t)=\frac{2\pi}{3},$$

\no such as shown in Fig.$4.18$ where each number on the side of a
vertex denotes $\rho_M(x)\mu(x)$ for $x=u,v,w$ and $t$.

\includegraphics[bb=5 5 100 120]{sgm72.eps}

\vskip 2mm

\c{\bf Fig.$4.18$}\vskip 2mm

\no We assume the linear assumption is holds in this planar map
geometry $(M,\mu)$. Then we get a classifications for $m$-points
in $(M,\mu)$ as follows.

$$V_{el}=\{\rm points \ in \ (uA\setminus\{A\})\bigcup
(uB\setminus\{B\})\bigcup (ut\setminus\{t\})\},$$

\no where $A$ and $B$ are euclidean points on $(u,w)$ and $(u,v)$,
respectively.

$$V_{eu}=\{A,B,t\}\bigcup(P\setminus E(K_4))$$

\no and

$$V_{hy}=\{\rm points \ in \ (wA\setminus \{A\})\bigcup(wt\setminus\{t\})
\bigcup wv\bigcup(tv\setminus\{t\})\bigcup(vB\setminus\{B\})\}.$$

Edges in $K_4$ are classified into $(u,t)\in C_E^1, (t,w),(t,v)\in
C_E^3$, $(u,w),(u,v)\in C_E^5$ and $(w,u)\in C_E^6$.

Various $m$-lines in this planar map geometry are shown in
Fig.$4.19$.

\includegraphics[bb=5 5 100 130]{sgm73.eps}

\vskip 2mm

\c{\bf Fig.$4.19$}\vskip 2mm

There are no $1$-polygons in this planar map geometry. One
$2$-polygon and various triangles are shown in Fig.$4.20$.

\includegraphics[bb=5 5 100 160]{sgm74.eps}

\vskip 2mm

\c{\bf Fig.$4.20$}

\vskip 4mm

\no{\bf Example $4.6.2$} \ {\it A wheel planar map $W_{1.4}$}

\vskip 3mm

We take a wheel $W_{1.4}$ embedded on a plane $\sum$ with vertices
$O$ and $u,v,w,t$ and angle factors

$$\mu(O)=\frac{\pi}{2},   \ {\rm and} \ \mu(u)=\mu(v)=\mu(w)=\mu(t)=\frac{4\pi}{3},$$

\no such as shown in Fig.$4.21$.

\includegraphics[bb=5 5 100 110]{sgm75.eps}

\vskip 2mm

\c{\bf Fig.$4.21$}

\vskip 2mm

There are no elliptic points in this planar map geometries.
Euclidean and hyperbolic points $V_{eu}, V_{hy}$ are as follows.

$$V_{eu}= P\bigcup\setminus(E(W_{1.4})\setminus\{O\})$$

\no and

$$V_{hy}= E(W_{1.4})\setminus\{O\}.$$

Edges are classified into $(O,u),(O,v),(O,w),(O,t)\in C_E^3$ and
$(u,v),(v,w),$ $(w,t),(t,u)\in C_E^6$. Various $m$-lines and one
$1$-polygon are shown in Fig.$4.22$ where each $m$-line will turn
to its opposite direction after it meeting $W_{1.4}$ such as those
$m$-lines $L_1, L_2$ and $L_4$, $L_5$ in Fig.$4.22$.

\includegraphics[bb=30 5 200 110]{sgm76.eps}

\vskip 2mm

\c{\bf Fig.$4.22$}

\vskip 4mm

\no{\bf Example $4.6.3$} \ {\it A parallel bundle in a planar map
geometry}

\vskip 3mm

We choose a planar ladder and define its angle factor as shown in
Fig.$4.23$ where each number on the side of a vertex $u$ denotes
the number $\rho_M(u)\mu(u)$. Then we find a parallel bundle
$\{L_i; 1\leq i\leq 6\}$ as those shown in Fig.$4.23$.

\includegraphics[bb=5 5 100 175]{sgm77.eps}

\vskip 2mm

\c{\bf Fig.$4.23$}

\vskip 5mm

\no{\bf \S $4.7$ \ Remarks and Open Problems}

\vskip 4mm

\no{\bf $4.7.1.$}  \ Unless the Einstein's relativity theory,
nearly all other branches of physics use an Euclid space as their
spacetime model. This has their own reason, also due to one's
sight because the moving of an object is more likely as it is
described by an Euclid geometry. As a generalization of an Euclid
geometry of plane by the Smarandache's notion, planar map
geometries were introduced in the references $[54]$ and $[62]$.
The same research can be also done for an Euclid geometry of a
space ${\bf R}^3$ and open problems are selected in the following.

\vskip 4mm

\no{\bf Problem $4.7.1$} \ {\it Establish Smarandache geometries
of a space ${\bf R}^3$ and classify their fundamental elements,
such as points, lines, polyhedrons, $\cdots$, etc..}

\vskip 4mm

\no{\bf Problem $4.7.2$} \ {\it Determine various surfaces in a
Smarandache geometry of a space ${\bf R}^3$, such as a sphere, a
surface of cylinder, circular cone, a torus, a double torus and a
projective plane, a Klein bottle, $\cdots$, also determine various
convex polyhedrons such as a tetrahedron, a pentahedron, a
hexahedron, $\cdots$, etc..}

\vskip 4mm

\no{\bf Problem $4.7.3$} \ {\it Define the conception of volume
and find formulae for volumes of convex polyhedrons in a
Smarandache geometry of a space ${\bf R}^3$, such as a
tetrahedron, a pentahedron or a hexahedron, $\cdots$, etc..}

\vskip 4mm

\no{\bf Problem $4.7.4$} \ {\it Apply Smarandache geometries of a
space ${\bf R}^3$ to find knots and characterize them.}

\vskip 4mm

\no{\bf $4.7.2.$} \ As those proved in Chapter $3$, we can also
research these map geometries on a locally orientable surfaces and
find its fundamental elements in a surface, such as a sphere, a
torus, a double torus, $\cdots$ and a projective plane, a Klein
bottle, $\cdots$, i.e., to establish an intrinsic geometry on a
surface. For this target, open problems for surfaces with small
genus should be solved in the first.

\vskip 4mm

\no{\bf Problem $4.7.5$} \ {\it Establish an intrinsic geometry by
map geometries on a sphere or a torus and find its fundamental
elements.}

\vskip 4mm

\no{\bf Problem $4.7.6$} \ {\it Establish an intrinsic geometry on
a projective or a Klein bottle and find its fundamental elements.}

\vskip 4mm

\no{\bf Problem $4.7.7$} \ {\it Define various measures of map
geometries on a locally orientable surface $S$ and apply them to
characterize the surface $S$.}

\vskip 4mm

\no{\bf Problem $4.7.8$} \ {\it Define the conception of curvature
for a map geometry $(M,\mu)$ on a locally orientable surface and
calculate the sum $\omega(M)$ of curvatures on all edges in $M$.}

\vskip 4mm

\no{\bf Conjecture $4.7.1$} \ {\it $\omega(M)=2\pi\chi(M)s(M)$,
where $s(M)$ denotes the sum of length of edges in $M$.}

\newpage

\no{\large\bf $5.$ Pseudo-Plane geometries}

\vskip 25mm

\no The essential idea in planar map geometries is to associate
each point in a planar map with an angle factor which turns
flatness of a plane to tortuous as we have seen in Chapter $4$.
When the order of a planar map tends to infinite and its diameter
of each face tends to zero (such planar maps exist, for example,
triangulations of a plane), we get a tortuous plane at the
limiting point, i.e., a plane equipped with a vector and straight
lines maybe not exist. We concentrate on discussing these
pseudo-planes in this chapter. A relation for integral curves with
differential equations is established, which enables us to find
good behaviors of plane curves.

\vskip 5mm

\no{\bf \S $5.1$ \ Pseudo-Planes}

\vskip 4mm

\no In the classical analytic geometry of plane, each point is
correspondent with the Descartes coordinate $(x,y)$, where $x$ and
$y$ are real numbers which ensures the flatness of a plane.
Motivated by the ideas in Chapters $3$ and $4$, we find a new kind
of planes, called {\it pseudo-planes} which distort the flatness
of a plane and can be applied to classical mathematics.

\vskip 4mm

\no{\bf Definition $5.1.1$}  \ {\it Let $\sum$ be an Euclid plane.
For $\forall u\in \sum$, if there is a continuous mapping $\omega:
u\rightarrow \omega(u)$ where $\omega(u)\in{\bf R}^n$ for an
integer $n, n\geq 1$ such that for any chosen number $\epsilon
>0$, there exists a number $\delta > 0$ and a
point $v\in \sum$, $\|u-v\|<\delta$ such that
$\|\omega(u)-\omega(v)\|< \epsilon$, then $\sum$ is called a
pseudo-plane, denoted by $(\sum,\omega)$, where $\|u-v\|$ denotes
the norm between points $u$ and $v$ in $\sum$.}

\vskip 3mm

An explanation for Definition $5.1.1$ is shown in Fig.$5.1$, in
where $n=1$ and $\omega(u)$ is an angle function $\forall u\in
\sum$.

\includegraphics[bb=10 10 100 120]{sgm78.eps}

\vskip 3mm

\c{\bf Fig.$5.1$}

\vskip 2mm

We can also explain $\omega(u)$, $u\in {\mathcal P}$ to be the
coordinate $z$ in $u=(x,y,z)\in{\bf R}^3$ by taking also $n=1$.
Thereby a pseudo-plane can be also seen as a projection of an
Euclid space ${\bf R}^{n+2}$ on an Euclid plane. This fact implies
that some characteristic of the geometry of space may reflected by
a pseudo-plane.

We only discuss the case of $n=1$ and explain $\omega(u), u\in
\sum$ being a periodic function in this chapter, i.e., for any
integer $k$, $4k\pi+\omega(u)\equiv\omega(u)(mod\ 4\pi)$. Not loss
of generality, we assume that $0< \omega(u)\leq 4\pi$ for $\forall
u\in \sum$. Similar to map geometries, points in a pseudo-plane
are classified into three classes, i.e., {\it elliptic points
$V_{el}$, euclidean points $V_{eu}$} and {\it hyperbolic points
$V_{hy}$}, defined by

$$V_{el} = \{u\in\sum| \omega(u) \ < 2\pi\},$$

$$V_{eu} = \{v\in\sum| \omega(v) \ =2\pi\}$$

\no and

$$V_{hy}= \{w\in\sum| \omega(w) \ >2\pi\}.$$

We define a sign function $sign(v)$ on a point of a pseudo-plane
$(\sum,\omega)$

\[
sign(v)=\left\{
\begin{array}{lr}
1, & {\rm if } \ v \ {\rm is \ elliptic},\\
0, & {\rm if } \ v \ {\rm is \ euclidean},\\
-1, & {\rm if } \ v \ {\rm is \ hyperbolic}.
\end{array}
\right.
\]

Then we get a result as in the following.

\vskip 4mm

\no{\bf Theorem $5.1.1$} \ {\it There is a straight line segment
$AB$ in a pseudo-plane $(\sum,\omega)$ if and only if for $\forall
u\in AB$, $\omega(u)=2\pi$, i.e., every point on $AB$ is
euclidean.}

\vskip 3mm

{\it Proof} \ Since $\omega(u)$ is an angle function for $\forall
u\in\sum$, we know that $AB$ is a straight line segment if and
only if for $\forall u\in AB$,

$$\frac{\omega(u)}{2}=\pi,$$

\no i.e., $\omega(u)=2\pi$, $u$ is an euclidean point.\quad\quad
$\natural$

Theorem $5.1.1$ implies that not every pseudo-plane has straight
line segments.

\vskip 4mm

\no{\bf Corollary $5.1.1$} \ {\it If there are only finite
euclidean points in a pseudo-plane $(\sum,\omega)$, then there are
no straight line segments in $(\sum,\omega)$.}

\vskip 4mm

\no{\bf Corollary $5.1.2$} \ {\it There are not always exist a
straight line between two given points $u$ and $v$ in a
pseudo-plane $({\mathcal P},\omega)$.}

\vskip 3mm

By the intermediate value theorem in calculus, we know the
following result for points in a pseudo-plane.

\vskip 4mm

\no{\bf Theorem $5.1.2$} \ {\it In a pseudo-plane $(\sum,\omega)$,
if $V_{el}\not=\emptyset$ and $V_{hy}\not=\emptyset$, then}

$$V_{eu}\not=\emptyset.$$

\vskip 3mm

{\it Proof} \ By these assumptions, we can choose points $u\in
V_{el}$ and $v\in V_{hy}$. Consider points on line segment $uv$ in
an Euclid plane $\sum$. Since $\omega(u)< 2\pi$ and $\omega(v)
> 2\pi$, there
exists at least a point $w,w\in uv$ such that $\omega(w)=2\pi$,
i.e., $w\in V_{eu}$ by the intermediate value theorem in calculus.
Whence, $V_{eu}\not=\emptyset$.\quad\quad $\natural$

\vskip 4mm

\no{\bf Corollary $5.1.3$} \ {\it In a pseudo-plane
$(\sum,\omega)$, if $V_{eu}=\emptyset$, then every point of
$(\sum,\omega)$ is elliptic or every point of $\sum$ is
hyperbolic.}

\vskip 3mm

According to Corollary $5.1.3$, pseudo-planes can be classified
into four classes as follows.

\vskip 3mm

{\bf $C_P^1$}(euclidean): {\it pseudo-planes whose each point is
euclidean.}

\vskip 2mm

{\bf $C_P^2$}(elliptic): {\it pseudo-planes whose each point is
elliptic.}

\vskip 2mm

{\bf $C_P^3$}(hyperbolic): {\it pseudo-planes whose each point is
hyperbolic.}

\vskip 2mm

{\bf $C_P^4$}(Smarandache's): {\it pseudo-planes in which there
are euclidean, elliptic and hyperbolic points simultaneously.}

\vskip 2mm

For the existence of an algebraic curve $C$ in a pseudo-plane
$(\sum,\omega)$, we get a criteria as in the following.

\vskip 4mm

\no{\bf Theorem $5.1.3$} \ {\it There is an algebraic curve
$F(x,y)=0$ passing through $(x_0,y_0)$ in a domain $D$ of a
pseudo-plane $(\sum,\omega)$ with Descartes coordinate system if
and only if $F(x_0,y_0)=0$ and for $\forall (x,y)\in D$,}

$$(\pi-\frac{\omega(x,y)}{2})(1+(\frac{dy}{dx})^2)=sign(x,y).$$

\vskip 3mm

{\it Proof} \ By the definition of pseudo-planes in the case of
that $\omega$ being an angle function and the geometrical meaning
of the differential value of a function at a point, we know that
an algebraic curve $F(x,y)=0$ exists in a domain $D$ of
$(\sum,\omega)$ if and only if

$$(\pi-\frac{\omega(x,y)}{2})=sign(x,y)\frac{d(\arctan(\frac{dy}{dx})}{dx},$$

\no for $\forall (x,y)\in D$, i.e.,

$$(\pi-\frac{\omega(x,y)}{2})=\frac{sign(x,y)}{1+(\frac{dy}{dx})^2},$$

\no such as shown in Fig.$5.2$, where $\theta=\pi-\angle 2+\angle
1$ ,$\lim\limits_{\triangle x\rightarrow 0}\theta= \omega(x,y)$
and $(x,y)$ is an elliptic point.

\includegraphics[bb=10 10 100 120]{sgm79.eps}

\vskip 3mm

\c{\bf Fig.$5.2$}

\vskip 2mm

\no Therefore we get that

$$(\pi-\frac{\omega(x,y)}{2})(1+(\frac{dy}{dx})^2)=sign(x,y).\quad\quad \natural$$

A plane curve $C$ is called {\it elliptic} or {\it hyperbolic} if
$sign(x,y)= 1$ or $-1$ for each point $(x,y)$ on $C$. We know a
result for the existence of an elliptic or a hyperbolic curve in a
pseudo-plane.

\vskip 4mm

\no{\bf Corollary $5.1.4$} \ {\it An elliptic curve $F(x,y)=0$
exists in a pseudo-plane $(\sum,\omega)$ with the Descartes
coordinate system passing through $(x_0,y_0)$ if and only if there
is a domain $D\subset\sum$ such that $F(x_0,y_0)=0$ and for
$\forall (x,y)\in D$,}

$$(\pi-\frac{\omega(x,y)}{2})(1+(\frac{dy}{dx})^2)=1$$

\no{\it and there exists a hyperbolic curve $H(x,y)=0$ in a
pseudo-plane $(\sum,\omega)$ with the Descartes coordinate system
passing through $(x_0,y_0)$ if and only if there is a domain
$U\subset\sum$ such that for $H(x_0,y_0)=0$ and $\forall (x,y)\in
U$,}

$$(\pi-\frac{\omega(x,y)}{2})(1+(\frac{dy}{dx})^2)=-1.$$

\vskip 3mm

Now construct a polar axis $(\rho,\theta)$ in a pseudo-plane
$(\sum,\omega)$. Then we get a result as in the following.

\vskip 4mm

\no{\bf Theorem $5.1.4$}\ {\it There is an algebraic curve
$f(\rho,\theta)=0$ passing through $(\rho_0,\theta_0)$ in a domain
$F$ of a pseudo-plane $(\sum,\omega)$ with a polar coordinate
system if and only if $f(\rho_0,\theta_0)=0$ and for $\forall
(\rho,\theta)\in F$,}

$$\pi-\frac{\omega(\rho,\theta)}{2}=sign(\rho,\theta)\frac{d\theta}{d\rho}.$$

\vskip 3mm

{\it Proof} \ Similar to the proof of Theorem $5.1.3$, we know
that $\lim\limits_{\triangle x\rightarrow 0}\theta= \omega(x,y)$
and $\theta=\pi-\angle 2+\angle 1$ if $(\rho,\theta)$ is elliptic,
or $\theta=\pi-\angle 1+\angle 2$ if $(\rho,\theta)$ is hyperbolic
in Fig.$5.2$. Whence, we get that

$$\pi-\frac{\omega(\rho,\theta)}{2}=sign(\rho,\theta)\frac{d\theta}{d\rho}.$$

\vskip 4mm

\no{\bf Corollary $5.1.5$} \ {\it An elliptic curve
$F(\rho,\theta)=0$ exists in a pseudo-plane $(\sum,\omega)$ with a
polar coordinate system passing through $(\rho_0,\theta_0)$ if and
only if there is a domain $F\subset\sum$ such that
$F(\rho_0,\theta_0)=0$ and for $\forall (\rho,\theta)\in F$,}

$$\pi-\frac{\omega(\rho,\theta)}{2}=\frac{d\theta}{d\rho}$$

\no{\it and there exists a hyperbolic curve $h(x,y)=0$ in a
pseudo-plane $(\sum,\omega)$ with a polar coordinate system
passing through $(\rho_0,\theta_0)$ if and only if there is a
domain $U\subset\sum$ such that $h(\rho_0,\theta_0)=0$ and for
$\forall (\rho,\theta)\in U$,}

$$\pi-\frac{\omega(\rho,\theta)}{2}=-\frac{d\theta}{d\rho}.$$

Now we discuss a kind of expressions in an Euclid plane ${\bf
R}^2$ for points in ${\bf R}^3$ and its characteristics.

\vskip 4mm

\no{\bf Definition $5.1.2$} \ {\it For a point $P=(x,y,z)\in {\bf
R}^3$ with center $O$, let $\vartheta$ be the angle of vector
$\overrightarrow{OP}$ with the plane $XOY$. Then define an angle
function $\omega:(x,y)\rightarrow 2(\pi-\vartheta)$, i.e., the
presentation of a point $(x,y,z)$ in ${\bf R}^3$ is a point
$(x,y)$ with $\omega(x,y)=2(\pi-\angle(\overrightarrow{OP},XOY))$
in a pseudo-plane $(\sum,\omega)$.}

\vskip 3mm

An explanation for Definition $5.2.1$ is shown in Fig.$5.3$ where
$\theta$ is an angle between the vector $\overrightarrow{OP}$ and
plane $XOY$.

\includegraphics[bb=10 10 100 120]{sgm80.eps}

\vskip 3mm

\c{\bf Fig.$5.3$}

\vskip 4mm

\no{\bf Theorem $5.1.5$} \ {\it Let $(\sum,\omega)$ be a
pseudo-plane and $P=(x,y,z)$ a point in ${\bf R}^3$. Then the
point $(x,y)$ is elliptic, euclidean or hyperbolic if and only if
$z>0$, $z=0$ or $z < 0$.}

\vskip 3mm

{\it Proof} \ By Definition $5.1.2$, we know that
$\omega(x,y)>2\pi$, $=2\pi$ or $< 2\pi$ if and only if $\theta
>0$,$= 0$ or $< 0$ since
$-\frac{\pi}{2}\leq\theta\leq\frac{\pi}{2}$. Those conditions are
equivalent to $z > 0$, $= 0$ or $< 0$. \quad\quad $\natural$

The following result reveals the shape of points with a constant
angle function value in a pseudo-plane $(\sum, \omega)$.

\vskip 4mm

\no{\bf Theorem $5.1.6$} \ {\it For a constant $\eta, 0 <\eta \leq
4\pi$, all points $(x,y,z)$ in ${\bf R}^3$ with $\omega(x,y)=\eta$
consist an infinite circular cone with vertex $O$ and an angle
$\pi-\frac{\eta}{2}$ between its generatrix and the plane $XOY$.}

\vskip 3mm

{\it Proof} \ Notice that $\omega(x_1,y_1)=\omega(x_2,y_2)$ for
two points $A, B$ in ${\bf R}^3$ with $A=(x_1,y_1,z_1)$ and
$B=(x_2,y_2,z_2)$ if and only if

$$\angle(\overrightarrow{OA},XOY)=\angle(\overrightarrow{OB},XOY)=\pi-\frac{\eta}{2},$$

\no that is, points $A$ and $B$ is on a circular cone with vertex
$O$ and an angle $\pi-\frac{\eta}{2}$ between
$\overrightarrow{OA}$ or $\overrightarrow{OB}$ and the plane
$XOY$. Since $z\rightarrow+\infty$, we get an infinite circular
cone in ${\bf R}^3$ with vertex $O$ and an angle
$\pi-\frac{\eta}{2}$ between its generatrix and the plane
$XOY$.\quad\quad $\natural$

\vskip 5mm

\no{\bf \S $5.2$ \ Integral Curves}

\vskip 4mm

\no An integral curve in an Euclid plane is defined by the next
definition.

\vskip 4mm

\no{\bf Definition $5.2.1$} \ {If the solution of a differential
equation}

$$\frac{dy}{dx}= f(x,y)$$

\no{\it with an initial condition $y(x_0)=y_0$ exists, then all
points $(x,y)$ consisted by their solutions of this initial
problem on an Euclid plane $\sum$ is called an integral curve.}

\vskip 3mm

By the ordinary differential equation theory, a well-known result
for the unique solution of an ordinary differential equation is
stated in the following. See also the reference $[3]$ for details.

\vskip 3mm

{\it If the following conditions hold:}\vskip 2mm

($i$) \ {\it $f(x,y)$ is continuous in a field $F$:}

$$F: x_0-a\leq x\leq x_0+a, \ \ y_0-b\leq y\leq y_0+b$$

($ii$) \ {\it there exist a constant $\varsigma$ such that for
$\forall (x,y), (x,\overline{y})\in F$,}

$$|f(x,y)-f(x,\overline{y})|\leq\varsigma|y-\overline{y}|,$$

\vskip 2mm

\no{\it then there is an unique solution}

$$y=\varphi(x), \ \ \ \varphi(x_0)=y_0$$

\no{\it for the differential equation}

$$\frac{dy}{dx}= f(x,y)$$

\no{\it with an initial condition $y(x_0)=y_0$ in the interval
$[x_0-h_0,x_0+h_0]$, where $h_0=\min(a,\frac{b}{M})$,
$M=\max\limits_{(x,y)\in R}|f(x,y)|$.}

\vskip 2mm

The conditions in this theorem are complex and can not be applied
conveniently. As we have seen in Section $5.1$ of this chapter, a
pseudo-plane $(\sum,\omega)$ is related with differential
equations in an Euclid plane $\sum$. Whence, by a geometrical
view, to find an integral curve in a pseudo-plane $(\sum,\omega)$
is equivalent to solve an initial problem for an ordinary
differential equation. Thereby we concentrate on to find integral
curves in a pseudo-plane in this section.

According to Theorem $5.1.3$, we get the following result.

\vskip 4mm

\no{\bf Theorem $5.2.1$} \ {\it A curve $C$,}

$$C=\{(x,y(x))| \frac{dy}{dx}=f(x,y), y(x_0)=y_0\}$$

\no{\it exists in a pseudo-plane $(\sum,\omega)$ if and only if
there is an interval $I=[x_0-h,x_0+h]$ and an angle function
$\omega:\sum\rightarrow{\bf R}$ such that}

$$\omega(x,y(x)) = 2(\pi - \frac{sign(x,y(x))}{1+f^2(x,y)})$$

\no{\it for $\forall x\in I$ with}

$$\omega(x_0,y(x_0))=2(\pi-\frac{sign(x,y(x))}{1+f^2(x_0,y(x_0))}).$$

\vskip 3mm

{\it Proof} \ According to Theorem $5.1.3$, a curve passing
through the point $(x_0,y(x_0))$ in a pseudo-plane $(\sum,\omega)$
if and only if $y(x_0)=y_0$ and for $\forall x\in I$,

$$(\pi-\frac{\omega(x,y(x))}{2})(1+(\frac{dy}{dx})^2)=sign(x,y(x)).$$

\no Solving $\omega(x,y(x))$ from this equation, we get that

$$\omega(x,y(x)) = 2(\pi - \frac{sign(x,y(x))}{1+(\frac{dy}{dx})^2})
= 2(\pi - \frac{sign(x,y(x))}{1+f^2(x,y)}).\quad\quad \natural$$

Now we consider curves with an constant angle function value at
each of its point.

\vskip 4mm

\no{\bf Theorem $5.2.2$} \ {\it Let $(\sum,\omega)$ be a
pseudo-plane. Then for a constant $0 <\theta \leq 4\pi$,}

($i$) \ {\it a curve $C$ passing through a point $(x_0,y_0)$ and
$\omega(x,y)=\eta$ for $\forall (x,y)\in C$ is closed without
self-intersections on $(\sum,\omega)$ if and only if there exists
a real number $s$ such that}

$$s\eta = 2(s-2)\pi.$$

($ii$) \ {\it a curve $C$ passing through a point $(x_0,y_0)$ with
$\omega(x,y)=\theta$ for $\forall (x,y)\in C$ is a circle on
$(\sum,\omega)$ if and only if}

$$\eta = 2\pi-\frac{2}{r},$$

\no{\it where $r=\sqrt{x_0^2+y_0^2}$, i.e., $C$ is a projection of
a section circle passing through a point $(x_0,y_0)$ on the plane
$XOY$.}

\vskip 3mm

{\it Proof} \ Similar to Theorem $4.3.1$, we know that a curve $C$
passing through a point $(x_0,y_0)$ in a pseudo-plane
$(\sum,\omega)$ is closed if and only if

$$\int\limits_{0}^s(\pi-\frac{\omega(s)}{2})ds=2\pi.$$

Now since $\omega(x,y)=\eta$ is constant for $\forall (x,y)\in C$,
we get that

$$\int\limits_{0}^s(\pi-\frac{\omega(s)}{2})ds=s(\pi-\frac{\eta}{2}).$$

\no Whence, we get that

$$s(\pi-\frac{\eta}{2})=2\pi,$$

\no i.e.,

$$s\eta = 2(s-2)\pi.$$

Now if $C$ is a circle passing through a point $(x_0,y_0)$ with
$\omega(x,y)=\theta$ for $\forall (x,y)\in C$, then by the Euclid
plane geometry we know that $s=2\pi r$, where
$r=\sqrt{x_0^2+y_0^2}$. Therefore, there must be that

$$\eta = 2\pi-\frac{2}{r}.$$

\no This completes the proof. \quad\quad $\natural$

Two spiral curves without self-intersections are shown in
Fig.$5.4$, in where $(a)$ is an input but $(b)$ an output curve.

\includegraphics[bb=10 10 100 140]{sgm81.eps}

\vskip 3mm

\c{\bf Fig.$5.4$}

\vskip 2mm

\no We call the curve in Fig.$5.4(a)$ an {\it elliptic in-spiral}
and Fig.$5.4(b)$ an {\it elliptic out-spiral}, correspondent to
the right hand rule. In a polar coordinate system $(\rho,
\theta)$, a spiral curve has equation

$$\rho =ce^{\theta t},$$

\no where $c,t$ are real numbers and $c >0$. If $t < 0$, then the
curve is an in-spiral as the curve shown in Fig.$5.4(a)$. If $t
>0$, then the curve is an out-spiral as shown in
Fig.$5.4(b)$.

For the case $t=0$, we get a circle $\rho = c$ (or $x^2+y^2=c^2$
in the Descartes coordinate system).

Now in a pseudo-plane, we can easily find conditions for in-spiral
or out-spiral curves. That is the following theorem.

\vskip 4mm

\no{\bf Theorem $5.2.3$} \ {\it Let $(\sum,\omega)$ be a
pseudo-plane and let $\eta, \zeta$ be constants. Then an elliptic
in-spiral curve $C$ with $\omega(x,y)=\eta$ for $\forall (x,y)\in
C$ exists in $(\sum,\omega)$ if and only if there exist numbers
$s_1>s_2>\cdots> s_l>\cdots$, $s_i >0$ for $i\geq 1$ such that}

$$s_i\eta < 2(s_i-2i)\pi$$

\no{\it for any integer $i, i\geq 1$ and an elliptic out-spiral
curve $C$ with $\omega(x,y)=\zeta$ for $\forall (x,y)\in C$ exists
in $(\sum,\omega)$ if and only if there exist numbers
$s_1>s_2>\cdots> s_l>\cdots$, $s_i >0$ for $i\geq 1$ such that}

$$s_i\zeta  > 2(s_i-2i)\pi$$

\no{\it for any integer $i, i\geq 1$.}

\vskip 3mm

{\it Proof} \ Let $L$ be an $m$-line like an elliptic in-spiral
shown in Fig.$5.5$, in where $x_1$, $x_2$,$\cdots$, $x_n$ are
non-euclidean points and $x_1x_6$ is an auxiliary line segment.

\includegraphics[bb=10 10 100 130]{sgm82.eps}

\vskip 3mm

\c{\bf Fig.$5.5$}

\vskip 2mm

\no Then we know that

$$\sum\limits_{i=1}^6(\pi-f(x_1)) < 2\pi,$$

$$\sum\limits_{i=1}^{12}(\pi-f(x_1)) < 4\pi,$$

$$\cdots\cdots\cdots\cdots\cdots\cdots.$$

Similarly from any initial point $O$ to a point $P$ far $s$ to $O$
on $C$, the sum of lost angles at $P$ is

$$\int\limits_{0}^s(\pi-\frac{\eta}{2})ds = (\pi-\frac{\eta}{2})s.$$

\no Whence, the curve $C$ is an elliptic in-spiral if and only if
there exist numbers $s_1>s_2>\cdots> s_l>\cdots$, $s_i >0$ for
$i\geq 1$ such that

$$(\pi-\frac{\eta}{2})s_1 < 2\pi,$$

$$(\pi-\frac{\eta}{2})s_2 < 4\pi,$$

$$(\pi-\frac{\eta}{2})s_3 < 6\pi,$$

$$\cdots\cdots\cdots\cdots\cdots\cdots,$$

$$(\pi-\frac{\eta}{2})s_l < 2l\pi.$$

\no Therefore, we get that

$$s_i\eta < 2(s_i-2i)\pi$$

\no for any integer $i, i\geq 1$.

Similarly, consider an $m$-line like an elliptic out-spiral with
$x_1$, $x_2$,$\cdots$, $x_n$ non-euclidean points. We can also
find that $C$ is an elliptic out-spiral if and only if there exist
numbers $s_1>s_2>\cdots> s_l>\cdots$, $s_i
>0$ for $i\geq 1$ such that

$$(\pi-\frac{\zeta}{2})s_1 > 2\pi,$$

$$(\pi-\frac{\zeta}{2})s_2 > 4\pi,$$

$$(\pi-\frac{\zeta}{2})s_3 > 6\pi,$$

$$\cdots\cdots\cdots\cdots\cdots\cdots,$$

$$(\pi-\frac{\zeta}{2})s_l > 2l\pi.$$

\no Whence, we get that

$$s_i\eta < 2(s_i-2i)\pi.$$

\no for any integer $i, i\geq 1$.\quad\quad $\natural$

Similar to elliptic in or out-spirals, we can also define a {\it
hyperbolic in-spiral} or {\it hyperbolic out-spiral} correspondent
to the left hand rule, which are mirrors of curves in Fig.$5.4$.
We get the following result for a hyperbolic in or out-spiral in a
pseudo-plane.

\vskip 4mm

\no{\bf Theorem $5.2.4$} \ {\it Let $(\sum,\omega)$ be a
pseudo-plane and let $\eta, \zeta$ be constants. Then a hyperbolic
in-spiral curve $C$ with $\omega(x,y)=\eta$ for $\forall (x,y)\in
C$ exists in $(\sum,\omega)$ if and only if there exist numbers
$s_1>s_2>\cdots> s_l>\cdots$, $s_i >0$ for $i\geq 1$ such that}

$$s_i\eta > 2(s_i-2i)\pi$$

\no{\it for any integer $i, i\geq 1$ and a hyperbolic out-spiral
curve $C$ with $\omega(x,y)=\zeta$ for $\forall (x,y)\in C$ exists
in $(\sum,\omega)$ if and only if there exist numbers
$s_1>s_2>\cdots> s_l>\cdots$, $s_i >0$ for $i\geq 1$ such that}

$$s_i\zeta  < 2(s_i-2i)\pi$$

\no{\it for any integer $i, i\geq 1$.}

\vskip 3mm

{\it Proof} \ The proof for $(i)$ and $(ii)$ is similar to the
proof of Theorem $5.2.3$. \quad\quad $\natural$

\vskip 5mm

\no{\bf \S $5.3$ \ Stability of a Differential Equation}

\vskip 4mm

\no  For an ordinary differential equation system

\begin{eqnarray*}
& \ & \frac{dx}{dt}=P(x,y),\\
& \ & \frac{dy}{dt}=Q(x,y),\ \ \ (A^*)
\end{eqnarray*}

\no where $t$ is a time parameter, the Euclid plane $XOY$ with the
Descartes coordinate system is called its a {\it phase plane} and
the orbit $(x(t),y(t))$ of its a solution $x=x(t), y=y(t)$ is
called an {\it orbit curve}. If there exists a point $(x_0,y_0)$
on $XOY$ such that

$$P(x_0,y_0)=Q(x_0,y_0)=0,$$

\no then there is an obit curve which is only a point $(x_0,y_0)$
on $XOY$. The point $(x_0,y_0)$ is called a {\it singular point of
$(A^*)$}. Singular points of an ordinary differential equation are
classified into four classes: {\it knot, saddle, focal} and {\it
central points}. Each of these classes are introduced in the
following.

\vskip 3mm

\no{\bf Class $1.$ Knots}

\vskip 2mm

A {\it knot} $O$ of a differential equation is shown in Fig.$5.6$
where $(a)$ denotes that $O$ is stable but $(b)$ is unstable.

\includegraphics[bb=10 10 100 130]{sgm83.eps}

\vskip 3mm

\c{\bf Fig.$5.6$}

\vskip 2mm

A {\it critical knot} $O$ of a differential equation is shown in
Fig.$5.7$ where $(a)$ denotes that $O$ is stable but $(b)$ is
unstable.

\includegraphics[bb=10 10 100 130]{sgm84.eps}

\vskip 3mm

\c{\bf Fig.$5.7$}

\vskip 2mm

A {\it degenerate knot} $O$ of a differential equation is shown in
Fig.$5.8$ where $(a)$ denotes that $O$ is stable but $(b)$ is
unstable.

\includegraphics[bb=10 10 100 130]{sgm85.eps}

\vskip 3mm

\c{\bf Fig.$5.8$}

\vskip 3mm

\no{\bf Class $2.$ Saddle points}

\vskip 2mm

A {\it saddle point} $O$ of a differential equation is shown in
Fig.$5.9$.

\includegraphics[bb=10 10 100 140]{sgm86.eps}

\vskip 3mm

\c{\bf Fig.$5.9$}

\vskip 3mm

\no{\bf Class $3.$ Focal points}

\vskip 2mm

A {\it focal point} $O$ of a differential equation is shown in
Fig.$5.10$ where $(a)$ denotes that $O$ is stable but $(b)$ is
unstable.

\includegraphics[bb=10 10 100 150]{sgm87.eps}

\vskip 3mm

\c{\bf Fig.$5.10$}

\vskip 3mm

\no{\bf Class $4.$ Central points}

\vskip 2mm

A {\it central point} $O$ of a differential equation is shown in
Fig.$5.11$, which is just the center of a circle.

\includegraphics[bb=10 10 100 120]{sgm88.eps}

\vskip 3mm

\c{\bf Fig.$5.11$}

\vskip 2mm

In a pseudo-plane $(\sum,\omega)$, not all kinds of singular
points exist. We get a result for singular points in a
pseudo-plane as in the following.

\vskip 4mm

\no{\bf Theorem $5.3.1$} \ {\it There are no saddle points and
stable knots in a pseudo-plane plane $(\sum,\omega)$.}

\vskip 3mm

{\it Proof} \ On a saddle point or a stable knot $O$, there are
two rays to $O$, seeing Fig.$5.6(a)$ and Fig.$5.10$ for details.
Notice that if this kind of orbit curves in Fig.$5.6(a)$ or
Fig.$5.10$ appears, then there must be that

$$\omega(O)=4\pi.$$

\no Now according to Theorem $5.1.1$, every point $u$ on those two
rays should be euclidean, i.e., $\omega(u)=2\pi$, unless the point
$O$. But then $\omega$ is not continuous at the point $O$, which
contradicts Definition $5.1.1$. \quad\quad $\natural$

If an ordinary differential equation system $(A^*)$ has a closed
orbit curve $C$ but all other orbit curves are not closed in a
neighborhood of $C$ nearly enough to $C$ and those orbits curve
tend to $C$ when $t\rightarrow+\infty$ or $t\rightarrow-\infty$,
then $C$ is called a {\it limiting ring of $(A^*)$} and {\it
stable} or {\it unstable} if $t\rightarrow+\infty$ or
$t\rightarrow-\infty$.

\vskip 4mm

\no{\bf Theorem $5.3.2$} \ {\it For two constants
$\rho_0,\theta_0$, $\rho_0 >0$ and $\theta_0\not= 0$, there is a
pseudo-plane $(\sum,\omega)$ with}

$$\omega(\rho,\theta)=2( \pi - \frac{\rho_0}{\theta_0\rho})$$

\no{\it or}

$$\omega(\rho,\theta)=2(\pi+ \frac{\rho_0}{\theta_0\rho})$$

\no{\it such that}

$$\rho = \rho_0$$

\no{\it is a limiting ring in $(\sum,\omega)$.}

\vskip 3mm

{\it Proof} \ Notice that for two given constants
$\rho_0,\theta_0$, $\rho_0 >0$ and $\theta_0\not=0$, the equation

$$\rho(t) =\rho_0e^{\theta_0\theta(t)}$$

\no has a stable or unstable limiting ring

$$\rho = \rho_0$$

\no if $\theta(t)\rightarrow 0$ when $t\rightarrow+\infty$ or
$t\rightarrow-\infty$. Whence, we know that

$$\theta(t)=\frac{1}{\theta_0}\ln\frac{\rho_0}{\rho(t)}.$$

\no Therefore,

$$\frac{d\theta}{d\rho}=\frac{\rho_0}{\theta_0\rho(t)}.$$

According to Theorem $5.1.4$, we get that

$$\omega(\rho,\theta)= 2(\pi - sign(\rho,\theta)\frac{d\theta}{d\rho}),$$

\no for any point $(\rho,\theta)\in\sum$, i.e.,

$$\omega(\rho,\theta)= 2(\pi - \frac{\rho_0}{\theta_0\rho})
$$

\no or

$$ \ \ \ \omega(\rho,\theta)=2(\pi+\frac{\rho_0}{\theta_0\rho}).  \ \ \
\natural$$

A general pseudo-space is discussed in the next section which
enables us to know the Finsler geometry is a particular case of
Smnarandache geometries.

\vskip 5mm

\no{\bf \S $5.4$ \ Remarks and Open Problems}

\vskip 4mm

\no Definition $5.1.1$ can be generalized as follows, which
enables us to enlarge our fields of mathematics for further
research.

\vskip 4mm

\no{\bf Definition $5.4.1$} \ {\it Let $U$ and $W$ be two metric
spaces with metric $\rho$, $W\subseteq U$. For $\forall u\in U$,
if there is a continuous mapping $\omega: u\rightarrow \omega(u)$,
where $\omega(u)\in {\bf R}^n$ for an integer $n, n\geq 1$ such
that for any number $\epsilon >0$, there exists a number $\delta
>0$ and a point $v\in W$, $\rho(u-v)<\delta$ such that
$\rho(\omega(u)-\omega(v))<\epsilon$, then $U$ is called a metric
pseudo-space if $U=W$ or a bounded metric pseudo-space if there is
a number $N > 0$ such that $\forall w\in W$, $\rho(w)\leq N$,
denoted by $(U,\omega)$ or $(U^-,\omega)$, respectively.}

\vskip 3mm

By choice different metric spaces $U$ and $W$ in this definition,
we can get various metric pseudo-spaces. For the case $n=1$, we
can also explain $\omega(u)$ being an angle function with $0
<\omega(u)\leq 4\pi$, i.e.,

\[
\omega(u)=\left\{
\begin{array}{ll}
\omega(u)(mod 4\pi), & {\rm if } \  {\rm u\in W},\\
2\pi, & {\rm if } \  {\rm u\in U\setminus W}\quad\quad (*)
\end{array}
\right.
\]

\no and get some interesting metric pseudo-spaces.

\vskip 4mm

\no{\bf $5.4.1.$ Bounded pseudo-plane geometries} \ Let $C$ be a
closed curve in an Euclid plane $\sum$ without self-intersections.
Then $C$ divides $\sum$ into two domains. One of them is finite.
Denote by $D_{fin}$ the finite one. Call $C$ a boundary of
$D_{fin}$. Now let $U=\sum$ and $W=D_{fin}$ in Definition $5.4.1$
for the case of $n=1$. For example, choose $C$ be a $6$-polygon
such as shown in Fig.$5.12$.

\includegraphics[bb=10 10 100 120]{sgm89.eps}

\vskip 3mm

\c{\bf Fig.$5.12$}

\vskip 2mm

\no Then we get a geometry $(\sum^-,\omega)$ partially euclidean
and partially non-euclidean.

\vskip 4mm

\no{\bf Problem $5.4.1$} \ {\it Similar to Theorem $4.5.2$, find
conditions for parallel bundles on $(\sum^-,\omega)$.}

\vskip 4mm

\no{\bf Problem $5.4.2$} \ {\it Find conditions for existing an
algebraic curve $F(x,y)=0$ on $(\sum^-,\omega)$.}

\vskip 4mm

\no{\bf Problem $5.4.3$} \ {\it Find conditions for existing an
integer curve $C$ on $(\sum^-,\omega)$.}

\vskip 4mm

\no{\bf $5.4.2.$ Pseudo-Space geometries} \ For any integer $m,
m\geq 3$ and a point $\overline{u}\in {\bf R}^m$. Choose $U=W={\bf
R}^m$ in Definition $5.4.1$ for the case of $n=1$ and
$\omega(\overline{u})$ an angle function. Then we get a
pseudo-space geometry $({\bf R}^m,\omega)$ on ${\bf R}^m$.

\vskip 4mm

\no{\bf Problem $5.4.4$} \ {\it Find conditions for existing an
algebraic surface $F(x_1,x_2,\cdots,x_m)=0$ in $({\bf
R}^m,\omega)$, particularly, for an algebraic surface
$F(x_1,x_2,x_3)=0$ existing in $({\bf R}^3,\omega)$.}

\vskip 4mm

\no{\bf Problem $5.4.5$} \ {\it Find conditions for existing an
integer surface in $({\bf R}^m,\omega)$.}

\vskip 3mm

If we take $U = {\bf R}^m$ and $W$ a bounded convex point set of
${\bf R}^m$ in Definition $5.4.1$. Then we get a bounded
pseudo-space $({\bf R}^{m-},\omega)$, which is partially euclidean
and partially non-euclidean.

\vskip 4mm

\no{\bf Problem $5.4.6$} \ {\it For a bounded pseudo-space $({\bf
R}^{m-},\omega)$, solve Problems $5.4.4$ and $5.4.5$ again.}

\vskip 4mm

\no{\bf $5.4.3.$ Pseudo-Surface geometries} \ For a locally
orientable surface $S$ and $\forall u\in S$, we choose $U=W=S$ in
Definition $5.4.1$ for $n=1$ and $\omega(u)$ an angle function.
Then we get a pseudo-surface geometry $(S,\omega)$ on the surface
$S$.

\vskip 4mm

\no{\bf Problem $5.4.7$} \ {\it Characterize curves on a surface
$S$ by choice angle function $\omega$. Whether can we classify
automorphisms on $S$ by applying pseudo-surface geometries
$(S,\omega)$?}

\vskip 3mm

Notice that Thurston had classified automorphisms of a surface $S,
\chi(S)\leq 0$ into three classes in $[86]$: {\it reducible,
periodic} or {\it pseudo-Anosov}.

If we take $U = S$ and $W$ a bounded simply connected domain of
$S$ in Definition $5.4.1$. Then we get a bounded pseudo-surface
$(S^{-},\omega)$.

\vskip 4mm

\no{\bf Problem $5.4.8$} \ {\it For a bounded pseudo-surface
$(S^{-},\omega)$, solve Problem $5.4.7$.}

\vskip 4mm

\no{\bf $5.4.4.$ Pseudo-Manifold geometries} \ For an $m$-manifold
$M^m$ and $\forall u\in M^m$, choose $U=W=M^m$ in Definition
$5.4.1$ for $n=1$ and $\omega(u)$ a smooth function. Then we get a
pseudo-manifold geometry $(M^m,\omega)$ on the $m$-manifold $M^m$.
This geometry includes the {\it Finsler geometry}, i.e., equipped
each $m$-manifold with a Minkowski norm defined in the following
($[13],[39]$).

A {\it Minkowski norm} on $M^m$ is a function $F:M^m\rightarrow
[0,+\infty)$ such that

($i$) \ $F$ is smooth on $M^m\setminus\{0\}$;

($ii$) \ $F$ is $1$-homogeneous, i.e.,
$F(\lambda\overline{u})=\lambda F(\overline{u})$ for
$\overline{u}\in M^m$ and $\lambda >0$;

($iii$) \ for $\forall y\in M^m\setminus\{0\}$, the symmetric
bilinear form $g_y: M^m\times M^m\rightarrow R$ with

$$g_y(\overline{u},\overline{v})=\frac{1}{2}
\frac{\partial^2F^2(y+s\overline{u}+t\overline{v})}{\partial
s\partial t}|_{t=s=0}$$.

\no is positive definite.

Then a {\it Finsler manifold} is a manifold $M^m$ and a function
$F:TM^m\rightarrow [0,+\infty)$ such that

($i$) \ $F$ is smooth on
$TM^m\setminus\{0\}=\bigcup\{T_{\overline{x}}M^m\setminus\{0\}:\overline{x}\in
M^m \}$;

($ii$) \ $F|_{T_{\overline{x}}M^m}\rightarrow[0,+\infty)$ is a
Minkowski norm for $\forall\overline{x}\in M^m$.\vskip 2mm

As a special case of pseudo-manifold geometries, we choose
$\omega(\overline{x})= F(\overline{x})$ for $\overline{x}\in M^m$,
then $(M^m,\omega)$ is a {\it Finsler manifold}, particularly, if
$\omega(\overline{x})=g_{\overline{x}}(y,y)= F^2(x,y)$, then
$(M^m,\omega)$ is a {\it Riemann manifold}. Thereby, {\it
Smarandache geometries, particularly pseudo-manifold geometries
include the Finsler geometry}.

Open problems for pseudo-manifold geometries are presented in the
following.

\vskip 4mm

\no{\bf Problem $5.4.9$} \ {\it Characterize these pseudo-manifold
geometries $(M^m,\omega)$ without boundary and apply them to
classical mathematics and to classical mechanics.}

\vskip 3mm

Similarly, if we take $U = M^m$ and $W$ a bounded submanifold of
$M^m$ in Definition $5.4.1$. Then we get a bounded pseudo-manifold
$(M^{m-},\omega)$.

\vskip 4mm

\no{\bf Problem $5.4.10$} \ {\it  Characterize these
pseudo-manifold geometries $(M^{m-},\omega)$ with boundary and
apply them to classical mathematics and to classical mechanics,
particularly, to hamiltonian mechanics.}

\end{document}